\documentclass{article}
\usepackage{epsfig}
\usepackage{latexsym}
\usepackage{amssymb}
\usepackage{amsmath}
\usepackage{amsfonts}
\usepackage{mathrsfs}
\usepackage[all]{xy}
\usepackage{texdraw}
\usepackage{graphicx}
\usepackage{psfrag}
\usepackage{pstricks}
\usepackage{makeidx}

\newtheorem{Theory}{Theory}[section] 
\newtheorem{theorem}[Theory]{Theorem}
\newtheorem{lemma}[Theory]{Lemma}
\newtheorem{technicalLemma}[Theory]{Technical Lemma}
\newtheorem{corollary}[Theory]{Corollary}
\newtheorem{proposition}[Theory]{Proposition}
\newtheorem{fact}{Fact}  
\newtheorem{remark}[Theory]{Remark} 
\newtheorem{question}{Question} 
\newtheorem{conjecture}[question]{Conjecture}

\newcommand{\Z}{\mathbb{Z}}
\newcommand{\N}{\mathbb{N}}
\newcommand{\Q}{\mathbb{Q}}
\newcommand{\R}{\mathbb{R}}
\newcommand{\be}{\begin{enumerate}}
\newcommand{\ee}{\end{enumerate}}
\newcommand{\bq}{\begin{question}}
\newcommand{\eq}{\end{question}}
\newcommand{\bcj}{\begin{conjecture}}
\newcommand{\ecj}{\end{conjecture}}
\newcommand{\bc}{\begin{corollary}}
\newcommand{\ec}{\end{corollary}}
\newcommand{\bl}{\begin{lemma}}
\newcommand{\el}{\end{lemma}}
\newcommand{\btl}{\begin{technicalLemma}}
\newcommand{\etl}{\end{technicalLemma}}
\newcommand{\bt}{\begin{theorem}}
\newcommand{\et}{\end{theorem}}
\newcommand{\bp}{\begin{proposition}}
\newcommand{\ep}{\end{proposition}}
\newcommand{\bft}{\begin{fact}}
\newcommand{\eft}{\end{fact}}
\newcommand{\brk}{\begin{remark}}
\newcommand{\erk}{\end{remark}}

\newcommand{\ploi}{\textrm{PL}_o(I)}

\newcommand{\supp}{\operatorname{Supp}}
\newcommand{\CS}{\mathfrak{C}}
\newcommand{\ut}[1]{\mathscr{T}_{#1}}
\newcommand{\leafl}{\lambda}

\title{Centralizers in the
R. Thompson group $V_n$}

\author{Collin Bleak \and Hannah Bowman
\and Alison Gordon \and Garrett Graham \and
Jacob Hughes \and Francesco Matucci \and
Eugenia Sapir} 

\begin{document}
\maketitle

\begin{abstract}
Let $n \ge 2$ and let $\alpha \in V_n$ be an element in the Higman-Thompson group $V_n$. 
We study the structure of the centralizer of $\alpha \in V_n$ 
through a careful analysis of the action of $\langle \alpha \rangle$ 
on the Cantor set $\CS$. 
We make use of revealing tree pairs as developed by Brin and Salazar from which we 
derive discrete train tracks and flow graphs to assist us in our analysis. A consequence of our structure 
theorem is that element centralizers are finitely generated. 
Along the way we give a short argument using 
revealing tree pairs which shows that cyclic groups are undistorted in $V_n$.
\end{abstract}

\section{Introduction}
\footnotetext[1]{Research supported 
by the National Science Foundation through the Research Experiences 
for Undergraduates grant at Cornell} 

In this paper, we produce a description of the structure of element
centralizers within the Higman--Thompson groups $V_n$. As a
corollary to our structure theorem we see that element centralizers in
the $V_n$ are finitely generated.  Finally, we give a separate short
argument which shows that all infinite-cyclic subgroups within $V_n$
are embedded without distortion.  (The groups $V_n$ are first
introduced by Higman in \cite{HigmanFPSG}, where $V_n = G_{n,1}$ in
his notation.  The R. Thompson group denoted $V$ in \cite{CFP} is
$V_2$ in our notation.)

For a given integer $n\geq 2$, our primary view of the group $V_n$ is as a
group of homeomorphisms acting on the Cantor set $\CS_n$ (seen as the
boundary of the rooted regular infinite $n$-ary tree $\mathscr{T}_n$).  This
point of view informs most of our work, where we use a close study of
the dynamics of subgroup actions on $\CS_n$ to derive our main result.

Centralizers in Thompson's group $F= F_2$ are first classified by
Guba and Sapir in \cite{GSDiagram} as a consequence of their
classification of element centralizers for diagram groups.  In related
but separate work, Brin and Squier in \cite{picric} describe roots and
centralizers in $\ploi$, the group of orientation-preserving,
piecewise linear homeomorphisms of the unit interval (from which it is
also easy to describe the element centralizers of $F_n$, although Brin
and Squier never formally do so).  Guba and Sapir also show in
\cite{GSFSubgroups} that element centralizers in $F$ are embedded
without distortion in $F$.  In Chapter 8 of the thesis \cite{MatucciThesis}, Bleak, Kassabov, and
Matucci classify centralizers in $T_2$ up to finite index.  This paper
can be seen as a continuation of the line of research leading to these
results.

  Higman's work in \cite{HigmanFPSG} also contains information about
  the structure of centralizers of elements in $V_n$ (see Theorem 9.9
  of \cite{HigmanFPSG}).  If one reads Higman's proof of Theorem 9.9
  carefully, one can derive with reasonable effort some of the
  information about the $\mathbb{Z}$ factors contained within the right hand direct product in \ref{BigTheorem}. 
  However, our own result contains significantly
  more detail about the overall structure of element centralizers than
  is contained within \cite{HigmanFPSG}.

  Separately, Tuna Altinel and Alexey Muranov (see
  \cite{AltinelMuranov}) use model theory to analyze aspects of the
  groups $F_n<T_n<V_n$.  In their work they compile some information
  with regard to element centralizers in these groups.  Their results
  with regards to element centralizers are very similar to what is
  known from the work of Brin and Squier in \cite{picric} and of
  Kassabov and Matucci in \cite{kmscpf}, and appears to be contained
  within the results of Higman from \cite{HigmanFPSG}.

  Mart\'inez-P\'erez and Nucinkis \cite{MartinezPerezNucinkis} recently studied generalizations
  of the groups of Higman, Thompson, Stein and Brin 
  and classified centralizers of finite subgroups
  in order to study finiteness properties of those groups via cohomology. 
  Their result generalizes
  the one in \cite{MatucciThesis} by Bleak, Kassabov, and
  Matucci and agrees with the one of the current paper when restricted
  to torsion elements.

The work in this paper uses in broad outline the approach of Bleak, Kassabov 
and Matucci to centralizers in Chapter 8 of \cite{MatucciThesis}, but we work in the more complex groups $V_n$, and thus
we need to employ a more complex set of tools in our analysis.  We
chose to use Brin's revealing pair technology (see \cite{brinHigherV})
for our supporting calculations (from amongst a fairly long list of
tools that provide similar data), and we developed discrete train tracks and our flow graph
objects to further support our intuitive understanding of how elements
of $V_n$ act on $\CS_n$.

One advantage of considering a revealing pair $(A,B,\sigma)$
representing an element $\alpha\in V_n$ over any random representative
tree pair is that it is easier to understand the orbit structure of
points in the Cantor set $\CS_n$ under the action of the group
$\langle \alpha\rangle$.  The paper \cite{salazarPaper} studies
revealing pairs in depth, and includes a solution of the conjugacy
problem for $V_n$ using revealing pairs (Salazar gives all arguments
in the context of $V_2$, but it is easy to see that her methods
extend to $n$-ary trees).  

We briefly describe some alternative technologies.  As alluded above,
one could use the seminormal and quasinormal forms from
\cite{HigmanFPSG} to gather much of the information we obtained from
revealing pairs. In fact, in response to early drafts of this article and conversations with
the first author of this article, Nathan Barker has replicated and extended 
many of the results
herein using Higman's technology and he has gone on
to work on the simultaneous conjugacy problem in $V_n$ \cite{BarkerThesis}. 
Another technology is the strand diagrams of Belk
and Matucci (see \cite{BelkMatucci}), which themselves are
refinements of Pride's pictures in \cite{Pr1, Pr2}. In turn, Pride's
pictures are essentially dual objects to the Dehn diagrams from
geometric group and semigroup theory (for instance, in this context,
one can study the related analysis of conjugacy in $F$ and other
diagram groups by Guba and Sapir in \cite{GSDiagram}).  In the end,
these tools all provide access to similar content.  We chose revealing
pairs as we were comfortable with calculations using them, and because
it was particularly easy to define our chief combinatorial objects, discrete train tracks and
flow graphs, from a revealing pair.

The dynamical information described by discrete train tracks and their 
corresponding flow graphs, forms a key ingredient in the proof by 
Bleak and Salazar \cite{bleakSalazar}
of the perhaps surprising result that $\mathbb{Z}^2*\mathbb{Z}$ does not embed in $V$. 
In particular, those authors make significant use of our flow graph technology in their analysis.

Sections 2 and 3, and subsection 4.1, serve as a mostly expository
introduction to calculations in the generalized R. Thompson groups
$F_n$, $T_n$, and $V_n$.  An informed reader in the area can likely
skip ahead to section 4.2, looking back to these sections on the rare
occasions in which a new term appears.

\bigskip
\noindent \textbf{Acknowledgments.}  The authors would like to thank
Claas R\"over and Martin Kassabov for interesting and helpful
conversations relating to this work.  We also owe a debt of thanks to
Mark Sapir for suggesting we look into the distortion of cyclic
subgroups in $V_n$. We thank Nathan Barker for interesting
conversations with respect to our flow graph object; it was by these
conversations that the authors were finally convinced of the need to
formally define flow-graphs, instead of just using them as informal
guides to intuition. Finally, the sixth author would like to thank
Conchita Mart\'inez-P\'erez for
pleasant conversations about her
work with Nucinkis \cite{MartinezPerezNucinkis}
where they discussed our independent approaches to the question, and the mutual confirmation of each other's results.

The authors also thank Cornell University for hosting, Robert
Strichartz for organizing, and the NSF for funding
the REU program at which most of this work took place.  We also thank
the Centre de Recerca Matem\`atica (CRM) in Barcelona, the University
of Virginia, and the University of Nebraska -- Lincoln for supporting
some of the authors and providing great conditions to finish this
work.

\subsection{A formal statement of results}

\newcommand{\ceq}{\mathrel{\mathop:}=}

If $\alpha\in V_n$, we set
\[
C_{V_n}(\alpha)\ceq\left\{\beta\in V_n \mid \alpha\beta=\beta\alpha\right\}.
\]
We call the group $C_{V_n}(\alpha)$ the \emph{centralizer in $V_n$ of $\alpha$}, as
is standard.

Our primary theorem is the following.

\begin{theorem} \label{BigTheorem}

Let $n>1$ be a positive integer and suppose $\alpha\in V_n$.  Then,
there are non-negative integers $s$, $t$, $m_i$, and $r_i$ and groups
$K_{m_i}$, $G_{n,r_i}$, $A_j$, and $P_t$, for $i\in \{1,\ldots, s\}$
and $j\in \{1,\ldots,t\}$, so that
\[
C_{V_n}(\alpha)\cong \left(\prod_{i = 1}^sK_{m_i}\rtimes G_{n,r_i} \right) 
\times
\left( \prod_{j= 1}^t \left(\left(A_j \rtimes \mathbb{Z}\right) \wr P_{q_j} \right) \right)
\]
\end{theorem}

We now explain the statement of this theorem in a bit more detail.

  The group $\langle\alpha\rangle$ acts on a subset of the nodes of
  the infinite $n$-ary tree.  The number $s$ represents the number of
  distinct lengths of finite cyclic orbits of nodes under this action.
  The value of $s$ is easy to compute from any given revealing pair
  representing $\alpha$.
  
For the action mentioned above, each $r_i$ is
  determined as a minimal number of nodes carrying a fundamental
  domain (for the action of a conjugate version of
  $\langle\alpha\rangle$) in the set of nodes supporting the cycles of
  length $m_i$.

For each $m_i$, we have $K_{m_i} = (Maps(\CS_n,\Z_{m_i}))^{r_i}$,
where $Maps(\CS_n,\Z_{m_i})$ is the group of continuous maps from
$\CS_n$ to $\Z_{m_i}$ under point-wise multiplication, and where
$\Z_{m_i}$ is the cyclic group $\Z/(m_i\Z)$ under the discrete
topology.  We note that $K_{m_i}$ is not finitely generated for
$m_i>1$.

The groups $G_{n,r_i}$ are the Higman-Thompson groups from
  \cite{HigmanFPSG}.

Given any element $\beta\in V_n$, one can associate an infinite
collection of finite labeled graphs (which we call \emph{flow
graphs}).  Flow graphs are labeled, directed, finite graphs and which
describe structural meta-data pertaining to the dynamics of certain subsets of
$\CS_n$ under the action of $\langle\alpha\rangle$. Flow graphs are themselves ``quotient
objects'' coming from \emph{discrete train tracks}, which are objects we introduce here
to model dynamics in the Cantor set much as regular train tracks model dynamics
in surface homeomorphism theory \cite{Thurston1, Thurston2, PennerHarer}

Components of a flow graph associated with $\alpha$ fall into equivalence
classes $\{ICC_i\}$ which model similar dynamics. The number $t$ is the number of equivalence classes of 
components carrying infinite orbits under the action of $\langle \alpha \rangle$. 
This number happens to be independent of the
representative flow graph chosen.

The right factor of the main direct product represents the restriction
of the centralizer of $\alpha$ to elements which are the identity away from
the closure of the region where $\langle \alpha\rangle$ acts with
non-finite orbits.

For each $j$, the supports of the elements of $ICC_j$ represent regions in the action of $\langle \alpha \rangle$
where the centralizers of $\alpha$, with restricted actions to these regions, are isomorphic.
For any such support of an element of $ICC_j$, the finite order elements of the restricted centralizer
form a group. We take $A_j$ to be a representative group from set of these isomorphic finite groups.

The group $P_{q_j}$ is the full symmetric group on the $q_j$ 
isomorphic flow graph components in $ICC_j$.

Recall now that by Higman in \cite{HigmanFPSG}, the groups $G_{k,r}$
are all finitely presented.  This fact, together with a short analysis
of the nature of the actions in the left-hand semi-direct products,
shows that each of the groups $K_{m_i}\rtimes G_{n,r_i}$ are finitely
generated (see Corollary \ref{finiteGenerationGmi} below).  As the
groups on the right-hand side of the central direct product are
manifestly finitely generated, we obtain the following corollary to
Theorem \ref{BigTheorem}.

\begin{corollary}\label{finiteGeneration}
Let $n>1$ be an integer and $\alpha\in V_n$.  The group
$C_{V_n}(\alpha)$ is finitely generated.
\end{corollary}

One could try to improve this last corollary to obtain a statement of
finite-presentation, which may be feasible.  Such a proof might be
accomplished through a careful study of presentations of the Higman
groups $G_{n,r}$.  In this direction, we ask the following.

\begin{question}
Must the group $C_{V_n}(\alpha)$ be finitely presented for every $\alpha\in V_n$?
\end{question}

We have also found some evidence supporting the possibility
that the answer to the following question is ``Yes.''

\begin{question}
Is it true that for each index $j$, the group $A_j$ is abelian?
\end{question}

It is not completely trivial to find an example element $\alpha\in V=
V_2$ where any $A_i$ is not cyclic.  In section
\ref{nontorsionSection} we give such an example where $t=1$ and $A_1\cong
\Z_2\times\Z_2$. 

A generator of each $\mathbb{Z}$ in the right hand terms
$A_i \rtimes \mathbb{Z}$ is given by a root of a restricted version of $\alpha$
which is restricted to act only on the support 
of an element in $ICC_j$.  We know that $A_i$ commutes with the restricted
version of $\alpha$ by definition, but it is not clear that $A_i$ will commute with
any valid choice of a generator for the $\mathbb{Z}$ term.

\begin{question}
Is it possible to replace the right hand terms $A_i \rtimes \mathbb{Z}$ with 
$A_i \times \mathbb{Z}$?
\end{question}

In the final section of the paper we prove the following theorem.
\begin{theorem}\label{distortion}
Suppose $\alpha\in V_n$ so that $\langle \alpha\rangle \cong \Z$.  The group $\langle \alpha\rangle$ is undistorted as a subgroup of $V_n$.
\end{theorem}

\subsection{Our general approach\label{bigApproach}}
We give a description of our approach to centralizers in the broadest
terms.

Let $\alpha\in V_n$.  Let $H = \langle \alpha \rangle$.  Define the
\emph{fundamental domain of the action of $H$} to be the space
$\CS_n/H$.  Note that this space is generally not very friendly,
i.e., it is typically not Hausdorff, but that fact will have almost no
bearing on our work.

If an element commutes with $\alpha$, it will induce
an action on $\CS_n/H$.

Thus, we get a short exact sequence:
\[
1\to \mathcal{K}\to C_{V_n}(\alpha)\to \mathcal{Q}\to 1
\]

Here, $\mathcal{K}$ represent the elements in $C_{V_n}(\alpha)$ which
act on $\CS_n$ in such a way that their induced action on $\CS_n/H$ is
trivial.  The group $\mathcal{Q}$ is the natural quotient of
$C_{V_n}(\alpha)$ by the image of the inclusion map $\mathcal{K}\to
C_{V_n}(\alpha)$.  Loosely speaking, elements of $\mathcal Q$ are
represented by elements of $V_n$ which act in the ``same'' way on each
``copy'' of the fundamental domain in $\CS_n$ (this of course is imprecise; there may be no embedding of a fundamental domain in
$\CS_n$).

The first author is indebted to Martin Kassabov for pointing out this
general structural approach to analyzing centralizers in groups of
homeomorphisms.

\section{Basic definitions}

Throughout this section, let us fix an integer $n>1$ for our
discussion.  We will also establish other conventions later that will
hold throughout the section, and not just in a particular
subsection.

We assume the reader is familiar with \cite{CFP}, and with the
definition of a tree pair representative of an element of $V = V_2$.
Nonetheless, we give an abbreviated tour through those definitions for
the non-experts (extending them to include the groups $F_n\le T_n\le V_n$) and
we state some essential lemmas which either occur in that source or
which are easily derived by the reader with full understanding of
these definitions.  We give some examples demonstrating most of that
language, and we add some new language to the lexicon mostly in
support of our own later definition of a flow graph.  In general the
reader experienced with R. Thompson group literature will find little
new material in this section and is encouraged to skip ahead,
returning only if he or she runs into an unfamiliar term in the later
parts of the paper.

\subsection{Trees and Cantor sets \label{TreesCantorSets}}
The only material in this section that may be unfamiliar to readers
conversant with R. Thompson group literature is some of the language
describing \emph{the Cantor set underlying a node of the tree
  $\ut{n}$} and related concepts.

Our primary perspective will be to consider $V_n$ as a group of
homeomorphisms of the Cantor set $\mathfrak{C}$.  In particular, $V_n$
should be thought to act as a group of homeomorphisms of the Cantor
set $\CS_n\cong \mathfrak{C}$.  That is, the version $\CS_n$ of the
Cantor set that is realized as the boundary of the standard infinite,
rooted $n$-ary tree $\ut{n}$.  While we assume the reader understands that
realization of the Cantor set, and also the terms ``node'' or
``vertex'', ``child'', ``parent'', ``ancestor'', ``descendant'',
``leaf'', and ``$n$-caret'' (or simply ``caret''), and similar related
language when referring to aspects of rooted $n$-ary trees, we give a short discussion below to establish some of our standard usage.

The left-to-right ordering of the children of a vertex in $\ut{n}$
allows us to give a name to each vertex in $\ut{n}$.  Given a vertex
$p$ in $\ut{n}$, there is a unique order preserving bijection
$\textrm{Ord}_p:Children(p)\to \left\{0,1,2,\ldots,n-1\right\}$.  Let
$v$ be a vertex of $\ut{n}$ and let $(v_i)_{i = 1}^m$ be the unique
descending path in $\ut{n}$ starting at the  root $r = v_1$ and ending at $v = v_m$,
then we name the vertex $v$ with the sequence
$(\textrm{Ord}_{v_i}(v_{i+1}))_{i = 1}^{m-1}$.  We will think of this
sequence as a string (ordered from left-to-right).

Given names of two vertices, we may concatenate these strings to produce the
name of a third vertex which will be a descendant of the first vertex

Below, we diagram an example of $\ut{2}$, with a finite tree $T$
highlighted within it.  The vertex $c$ is
a leaf of $T$, and in both $T$ and $\ut{2}$, $c$ is a child of
$b$ which is a child of $a$.  The vertex $a$ is an ancestor of $c$ and
$c$ is a descendant of $a$.  The name of the vertex labeled by $c$ is
$010$.

\begin{center}
\psfrag{1}[c]{$1$}
\psfrag{0}[c]{$0$}
\psfrag{a}[c]{$a$}
\psfrag{b}[c]{$b$}
\psfrag{c}[c]{$c$}
\psfrag{T}[c]{$T$}
\includegraphics[height=200pt,width=350 pt]{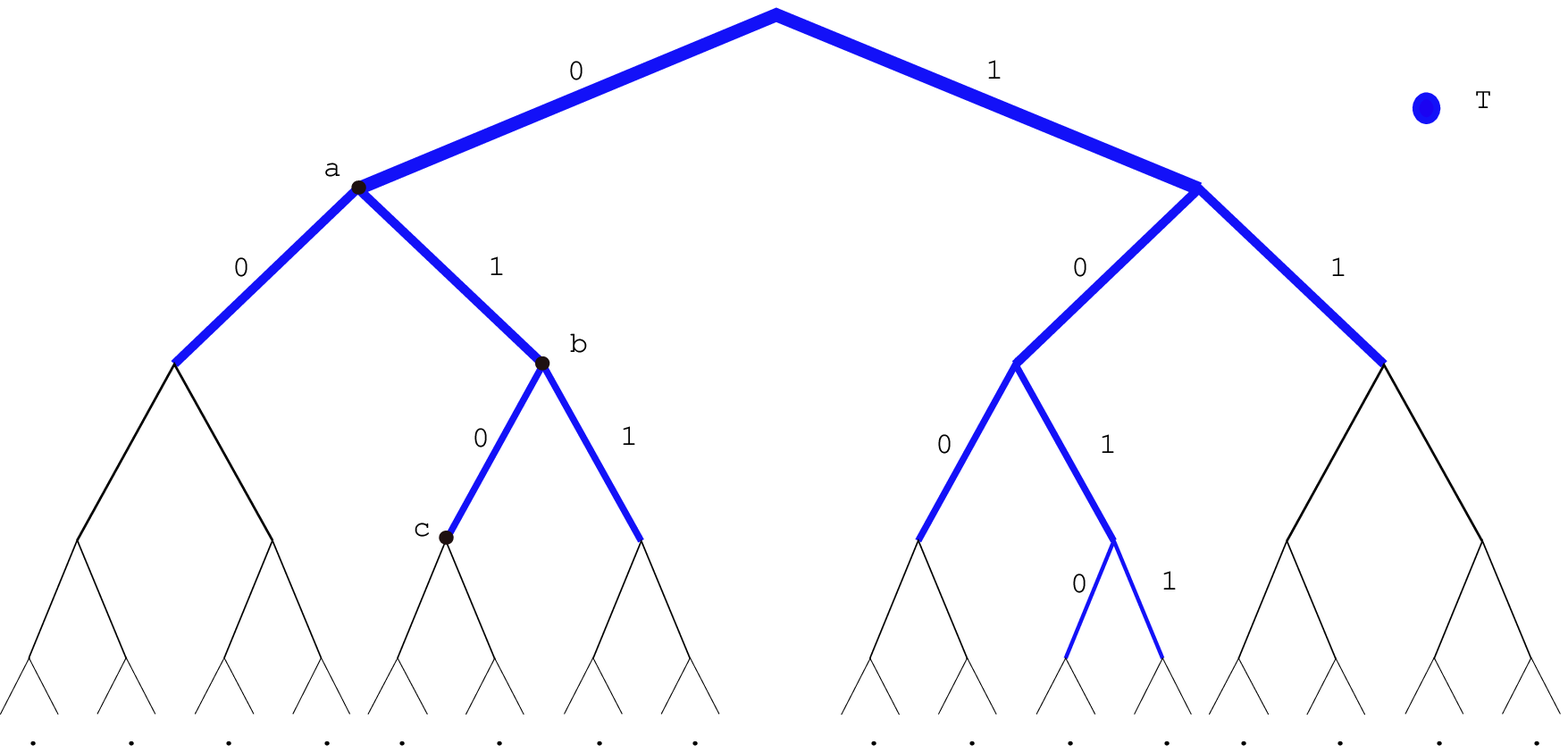}
\end{center}

We can view a finite rooted $n$-ary tree $T$ as an instruction on how
to partition the Cantor set $\CS_n$.  Consider the natural embedding
of $T$ into the tree $\mathscr{T}_n$, where we send the root of $T$ to
the root of $\mathscr{T}_n$, and we preserve orders of children.  For
instance, as in the diagram above.  Now, consider the set $P_n$ of all
infinite descending paths in $\mathscr{T}_n$ which start at the root
of $\mathscr{T}_n$.  If we consider each ordered $n$-caret in
$\mathscr{T}_n$ as an instruction to pass through another inductive
subdivision of the unit interval in the formation process of $\CS_n$,
then each element in $P_n$ can be thought of as limiting on an element
of $\CS_n$.  We thus identify $P_n$ with $\CS_n$. The set $P_n$ will
now be considered to be topologized using the induced topology from
the metric space topology of the unit interval.  Now if we consider a
vertex $c$ of $T$, we can associate $c$ with the subset of $\CS_n$
corresponding to the paths in $P_n$ which pass through $c$.  We will
call this \emph{the Cantor set under $c$}, and we will call any such
subset of $\CS_n$ an \emph{interval of $\CS_n$}.  It is immediate that
any interval in $\CS_n$ is actually homeomorphic with $\CS_n$.  Given a node $c$, the Cantor set under $c$ is also commonly called a cone neighborhood in $\CS_n$, and by definitions these sets form the standard basis for the product topology on $P_n = \{0,1,\ldots,n-1\}^{\omega} \cong \CS_n$. 
Thus, the
leaves of $T$ partition the set $P_n$, and they also partition $\CS_n$, into a set of open basis sets.
We will call this \emph{the partition of $\CS_n$ associated with the
  tree $T$}.

Extending the language of the previous paragraph, given a set $X \subset \CS_n$,
we will call any node $c$ of the universal tree $\mathscr{T}_n$ which has its underlying
set contained in $X$, a \emph{node of $X$}.

Of course, for all integers $m,n >1$, we have that $\CS_m\cong
\CS_n\cong \CS_2 = \CS$.

Using the example above, the interval of $\CS_2$ under $c$ is
$\CS_2\cap [2/9,7/27]$.  In discussion, we will generally not
distinguish between a vertex of $\ut{n}$ and the interval under it.

\begin{remark}
Any finite union of disjoint intervals in $\CS_n$ is homeomorphic with $\CS_n$.
\end{remark}

\subsection{Elements of $V_n$, $T_n$, and $F_n$
\label{sec:examples-V_n-T_n-F_n}}
Some of the language in this subsection is unusual, although the
general content will be familiar to all readers with knowledge of the
R. Thompson groups.

An element of $\textrm{Homeo}(\CS_n)$ is \emph{allowable} if it can be
represented by an allowable triple $(A, B, \sigma)$.  The triple
$(A,B,\sigma)$ is \emph{allowable} if there is a positive integer $m$
so that $A$ and $B$ are rooted, finite, $n$-ary trees with the same
number $m$ of leaves, and $\sigma\in\Sigma_m$, the permutation group
on the set $\left\{1,2,3,\ldots,m\right\}$.  We explain below how to
build the homeomorphism $\alpha$ which is represented by such an
allowable triple $(A,B,\sigma)$.  We then call $(A,B,\sigma)$ a
\emph{representative tree pair for $\alpha$}.  The group $V_n$
consists of the set of all allowable homeomorphisms of $\CS_n$ under
the operation of composition.

We are now ready to explain how an allowable triple $(A,B,\sigma)$
defines an allowable homeomorphism $\alpha$ of $\CS_n$.  Suppose $A$ and
$B$ both have $m$ leaves, for some integer $m>1$.  We consider $A$ to
represent the domain of $\alpha$, and $B$ to represent the range of $\alpha$.
We take the leaves of $A$ and $B$ and number them in their natural
left-to-right ordering from $1$ to $m$.  For each index $i$ with
$1\leq i\leq m$, we map the interval under leaf $i$ of $A$ to the
interval under the leaf $i\sigma$ of $B$ using an
orientation-preserving affine homeomorphism (that is, the
homeomorphism of the two Cantor sets underlying these leaves defined
by a restriction and co-restriction of an affine map with positive
slope from the real numbers $\R$ to $\R$). Note for $\CS_n$ the slopes of such maps will be
integral powers of $2n-1$. We will say that \emph{the
leaf $i$ of $A$ is mapped to the leaf $i\sigma$ of $B$ by $\alpha$}.

The next remark follows immediately from the discussion above:
\begin{remark}\label{neutralTravelers}
Suppose $v\in V_n$, and $v$ is represented by a tree pair $(A,B,\sigma)$.  If $\eta$ is a node in $\mathscr{T}_n$ so that $\eta$ is either a leaf of $A$ or a descendant node in $\mathscr{T}_n$ of a leaf of $A$, then $v$ will carry the Cantor set underlying $\eta$ affinely and bijectively in an order-preserving fashion to the Cantor set underlying some node $\tau$ of $\mathscr{T}_n$, where $\tau$ is either a leaf of $B$ or a descendant of a leaf of $B$.
\end{remark}

The diagram below illustrates an example for $V_2=V$.  The tree $A$ is
on the left, and the tree $B$ is on the right.  Note that we have
re-decorated the leaves of $B$ with the numbering from $\sigma$.  The
diagram beneath the tree pair is intended to indicate where intervals
of $\CS_2$ are getting mapped, where the intervals $D$ represent the
domain and the intervals $R$ represent the range.  The lower diagram
is superfluous when defining a general element of $V$.

\begin{center}
\psfrag{1}[c]{$1$}
\psfrag{2}[c]{$\,\,2$}
\psfrag{3}[c]{$3$}
\psfrag{4}[c]{$4$}
\psfrag{5}[c]{$5$}
\psfrag{6}[c]{$6$}
\psfrag{D}[c]{$D$}
\psfrag{R}[c]{$R$}
\includegraphics[height=200pt,width=300 pt]{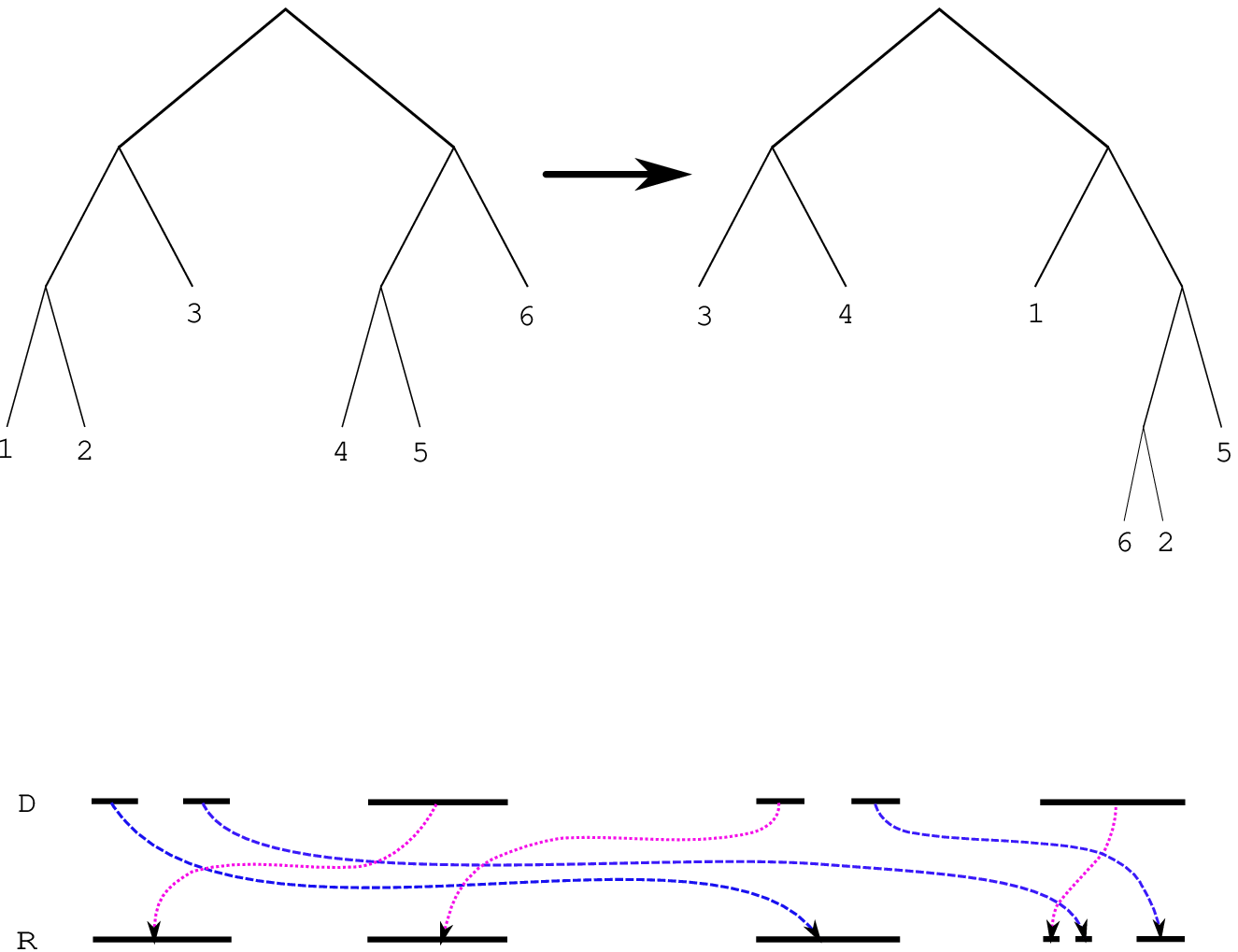}
\end{center}

In general, as we can decorate our tree leaves with numbers to
indicate the bijection, \emph{we will now re-define the phrase ``tree pair''
throughout the remainder of the paper to mean an allowable triple}.  We
will still discuss the permutation of a tree pair as needed.

As mentioned in the introduction, there are groups $F_n\leq T_n\leq
V_n$. If we only allow cyclic permutations, we get the group $T_n$.
If our permutation is trivial, we get $F_n$.  Thus, we can think of
$F_n$ as a group of piecewise-linear homeomorphisms of the interval
$[0,1]$, while $T_n$ can be thought of as a group of piecewise-linear
homeomorphisms of the circle $S^1$.

Here is an example element $\theta \in T = T_2$, which we will be considering again later.

\begin{center}
\psfrag{1}[c]{$1$}
\psfrag{2}[c]{$2$}
\psfrag{3}[c]{$3$}
\psfrag{4}[c]{$4$}
\psfrag{5}[c]{$5$}
\psfrag{6}[c]{$6$}
\psfrag{7}[c]{$7$}
\psfrag{D}[c]{$D$}
\psfrag{R}[c]{$R$}
\includegraphics[height=200pt,width=300 pt]{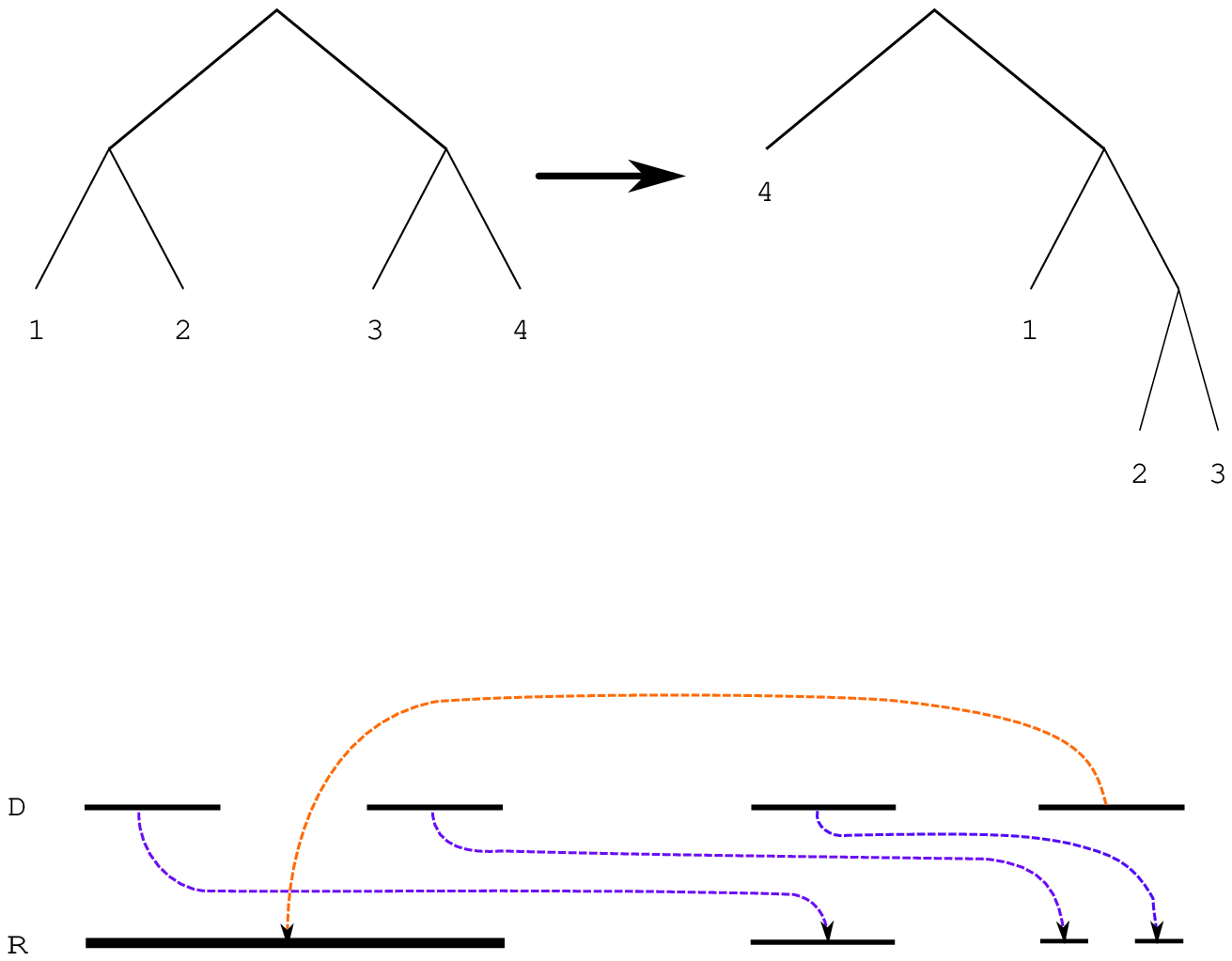}
\end{center}

Below is an example element from $F=F_2$.

\begin{center}
\psfrag{1}[c]{$1$}
\psfrag{2}[c]{$2$}
\psfrag{3}[c]{$3$}
\psfrag{4}[c]{$4$}
\psfrag{5}[c]{$5$}
\psfrag{6}[c]{$6$}
\psfrag{7}[c]{$7$}
\psfrag{D}[c]{$D$}
\psfrag{R}[c]{$R$}
\includegraphics[height=200pt,width=300 pt]{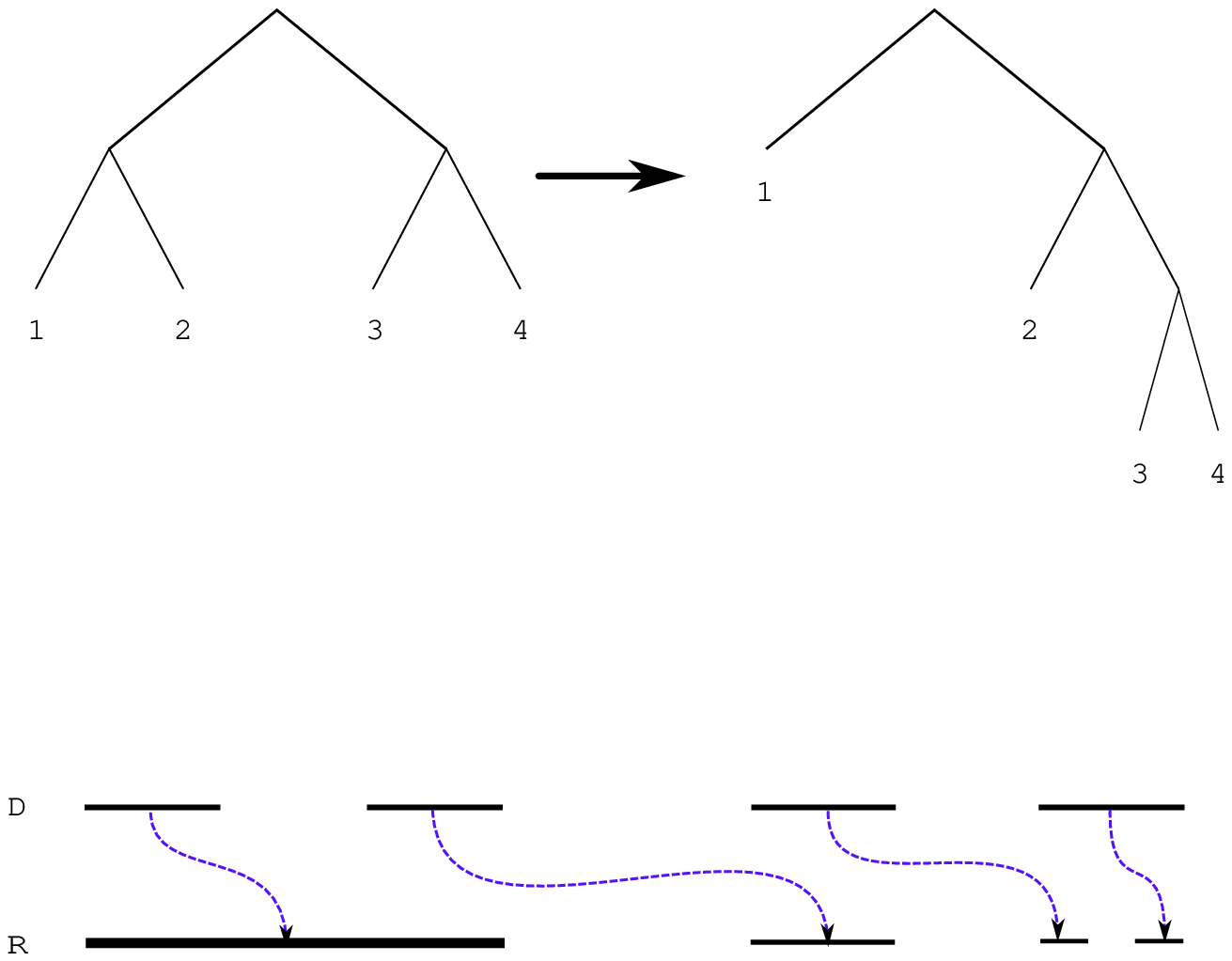}
\end{center}

\subsection{Multiplication in $V_n$}
There is nothing new in this subsection; readers familiar
with the R. Thompson groups should skip ahead.

Naturally, the group operation of $V_n$ is given by composition of
functions.  Building compositions using tree pairs is not hard.  The
process is enabled by the fact that there are many representative tree
pairs for an element of $V_n$.

We give an example of what we mean by the last sentence.  Given a
representative tree pair $(A,B,\sigma)$, one can use a leaf $i$ of $A$
to be a root of an extra $n$-caret, creating $A'$, and one can build a
tree $B'$ from $B$ by replacing the leaf $i\sigma$ with an $n$-caret.
Label the leaves of $A'$ in increasing order, and the leaves of $B'$
using the induced labeling from the permutation $\sigma$ on all
leaves of $B'$ that are also leaves of $B$.  For the other leaves of
$B'$, use the labeling, in order, of the $n$-caret in $A'$ that is
not a $n$-caret of $A$.  Let us call the permutation we have built
from the leaves of $A'$ to the leaves of $B'$ by $\sigma'$.  The
process of replacing the tree pair $(A,B,\sigma)$ by $(A',B',\sigma')$
is called a \emph{simple augmentation}.  

If the reverse process can be carried out (that is, deleting an
$n$-caret from both the domain and range trees and re-labeling the
permutation so that our initial two trees appear as a simple
augmentation of our resulting tree pair) then we call this process a
\emph{simple reduction}.  If we carry out either a simple augmentation
or a simple reduction, we may instead say we have done a \emph{simple
modification} to our initial tree pair.

The following lemma is a straightforward consequence of the standard
fact that any element in $V_n$ has a unique tree pair representation
that will not admit any simple reductions.

\begin{lemma}

Any two representative tree pairs of a particular element of $V_n$ are
connected by a finite sequence of simple modifications.  

\end{lemma}

We are now ready to carry out multiplication of tree pairs.  The
essence of the idea is to augment the range tree of the first element
and the domain tree of the second element until they are the same tree
(and carry out the necessary augmentations throughout both tree pairs),
and then re-label the permutations, using the labeling of the range
tree of the first element to seed the re-labeling of the permutation
in the second element's tree pair. At this junction, the two inner
trees are completely identical, and can be removed.  The diagram
below demonstrates this process.

\begin{center}
\psfrag{1}[c]{$1$}
\psfrag{2}[c]{$2$}
\psfrag{3}[c]{$3$}
\psfrag{4}[c]{$4$}
\psfrag{5}[c]{$5$}
\psfrag{6}[c]{$6$}
\includegraphics[height=350pt,width=250 pt]{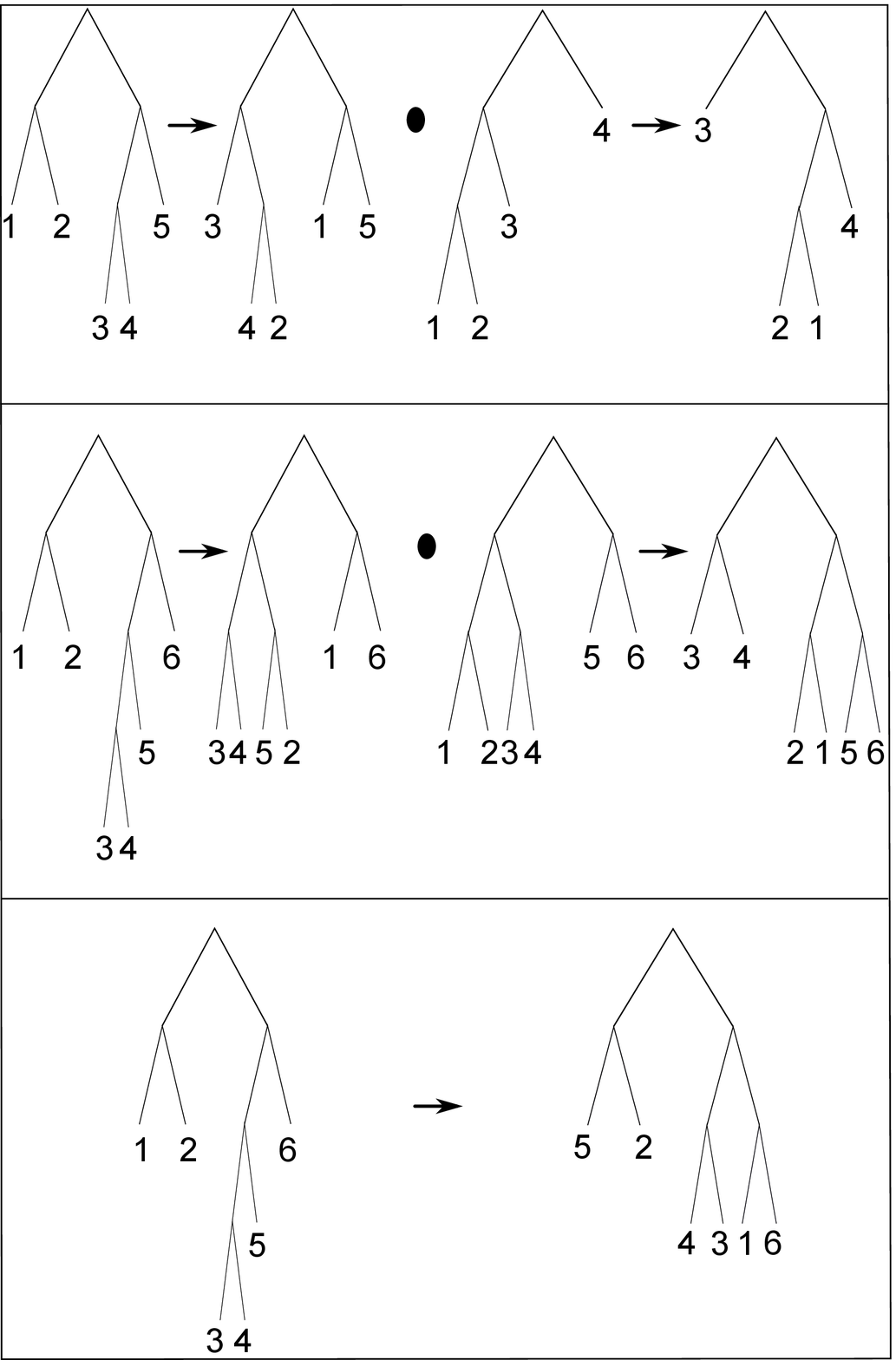}
\end{center}

Thus, we now know what $V_n$ is, and we can represent and
multiply its elements.

\section{Conjugation, roots, and centralizers}
Recall that we will use right action notation to describe how $V_n$
  acts on $\mathfrak{C}=\CS_n$, as below.

The following discussion is completely basic, and just carries out
some straightforward points from permutation group theory.

  Let $c\in\CS_n$ and $\alpha$, $\beta\in V_n$.  In particular,
  $\alpha:\CS_n\to\CS_n$.  

  \begin{itemize}

  \item We denote the image of $c$ under $\alpha$
    by $c\alpha$.
  \item For conjugation, we denote the conjugate of $\alpha$ by
  $\beta$ by $\alpha^\beta$, and note therefore that $c\alpha^{\beta}
  = c\beta^{-1}\alpha\beta$.
    
  \end{itemize}

  We will often need to discuss what is moving under an
  action, so we need a definition as well.  We define the support of
  an element $\alpha\in V_n$, denoted by $\supp{\alpha}$, as below.

\[
\supp(\alpha) = \left\{c\in \CS_n \mid c\alpha \neq
c\right\}.
\]

Note that this is distinct from the standard analysis version of
support, which would extend the definition to include the closure of
the set of points which are moving.

The following lemma is now standard from the theory of permutation groups.
\begin{lemma}

\[
\supp(\alpha^\beta)=\supp(\alpha)\beta
\]

\end{lemma}

Recall that given $\alpha\in V_n$ we have 
 
\[
C_{V_n}(\alpha)=\left\{\beta\in V_n \mid \alpha\beta = \beta\alpha\right\}.
\]

We point out the following obvious facts.

\begin{remark}
\label{conjugationRemark}
Suppose $\alpha$, $\beta$, and $\gamma\in V_n$.
\begin{enumerate}
\item If $\alpha\beta =\beta\alpha$, then $\alpha^{\beta} =
\beta^{-1}\alpha\beta = \alpha$.
\item If $\beta^k = \alpha$ for some integer $k$, then $\beta\in C_{V_n}(\alpha)$.
\item We have $C_{V_n}(\alpha)\cong (C_{V_n}(\alpha))^{\gamma}=
C_{V_n}(\alpha^{\gamma})$.
\end{enumerate}
\end{remark}

We will use the third point above repeatedly to replace an element whose
centralizer we are studying, by a conjugate element which admits
a simpler tree-pair representative (simplifying our
analysis without affecting centralizer structure).

\section{Revealing pairs and related objects\label{revealingPairs}}
Subsection \ref{rev-pair-def} should be considered as mostly expository; it will
contain definitions and lemmas from within \cite{brinHigherV} and
\cite{salazarPaper}.  We give very detailed examples of all of the
concepts therein which are of use in our context.  The following subsections on
discrete train tracks, laminations and flow graphs, on the other hand, are entirely new.

Many tree pairs exist to represent a single element.  Some tree pairs are more useful than others when it comes to discerning aspects of the dynamics of the element's action on the Cantor set.  Consider the element $\theta \in T$
which was defined earlier as the second tree pair diagram in Subsection \ref{sec:examples-V_n-T_n-F_n}.

Thought of as a homeomorphism of the circle, $\theta$ has a rotation
number, which, roughly stated, measures the average rotation of the
circle under the action of $\theta$.  (Rotation numbers are a
beautiful idea of Poincar\'{e}, and they are extremely useful in the
analysis of circle maps).  It is known, initially by work of Ghys and
Sergiescu (see \cite{GhysSergiescu}), that the rotation number of any
element of $T$ is rational.  Thus, $\theta$ has a rational rotation
number $p/q$ (in lowest terms).  A lemma of Poincar\'{e} now shows
that some point on the circle will have a periodic orbit, with period
$q$.  If the reader examines the tree pair representative for
$\theta$, doubtless he or she will discover a point that travels on a
finite orbit of length five.  In fact, all points of the circle travel
on their own periodic orbit of length five, as $t$ is torsion with
order five.  Under any reasonable definition of rotation number,
$\theta$ must have rotation number $2/5$, since the reader can observe
that after five iterations of $\theta$, the circle will have rotated a
total of two whole times around.

The previous example is intended to point to the fact that some tree
pairs somehow hide the dynamics of an element's action.

In this section, we will describe revealing pairs, which are tree
pairs that easily yield up all information about the orbit dynamics in
$\CS_n$ under the action of a cyclic subgroup of $V_n$ (or in the
circle or the interval, in the cases of cyclic subgroups of $T_n$ or
$F_n$, respectively).

\subsection{Revealing pair definitions \label{rev-pair-def}}
Throughout this section, we will work with some nontrivial $\alpha\in V_n$.
We will assume that the tree pair $(A,B,\sigma)$ represents $\alpha$.

Consider $A$ and $B$ as finite rooted
$n$-ary sub-trees (with roots at the root of $\mathscr{T}_n$) of
$\mathscr{T}_n$.  

We will call the set of vertices of $\mathscr{T}_n$ which are leaves
of both $A$ and $B$ the \emph{neutral leaves of $(A,B,\sigma)$}.  We
will simply call these the neutral leaves, if the tree pair is
understood.
 
As both $A$ and $B$ have a root and neither are empty, we can
immediately form the tree $C = A\cap B$.  It is immediate that the
neutral leaves are leaves of $C$, but if $A\neq B$, then $C$ will have
other leaves as well.

We can make the sets $X = A-B$ and $Y = B-A$.  The closures
$\overline{X}$ and $\overline{Y}$ in $\mathscr{T}_n$ are both finite
disjoint unions rooted $n$-ary trees (their roots are not sitting at
the root of $\mathscr{T}_n$), where the number of carets in
$\overline{X}$ is the same as the number of carets of $\overline{Y}$
(this number could be zero if $A=B$).  We call $\overline{X}$ and
$\overline{Y}$ \emph{a difference of carets} for $A$ and $B$.

In the remainder, when we write $D-E$, where $D$ and $E$ are trees, we
actually want to take the difference of carets, so that our result
will be a collection of rooted trees.

While $\langle \alpha\rangle$ acts on the Cantor set, it also induces a
``partial action'' on an infinite subset of the vertices of
$\mathscr{T}_n$, as we explain in this paragraph.  Since the interval
under a leaf $\lambda$ of $A$ is taken affinely to the interval under
a leaf of $B$ (we will denote this leaf by $\lambda \alpha$), we see that
$\alpha$ induces a map from the vertices of $\ut{n}$ under $\leafl$
to the vertices of $\ut{n}$ under $\leafl \alpha$.  In particular, the full
sub-tree in $\ut{n}$ with root $\leafl$ is taken to the full sub-tree in
$\ut{n}$ with root $\leafl \alpha$ in order preserving fashion.  We note in
passing that we cannot extend this to a true action on the vertices of
$\ut{n}$; if we consider a vertex $\eta$ in $\ut{n}$ which is above a leaf
of $A$, the map $\alpha$ may take the interval underlying $\eta$ and map it
in a non-affine fashion across multiple intervals in $\CS_n$.

As an example of the behavior mentioned above, consider the element
$\theta$ again.  The parent vertex of the domain leaves labeled 3 and 4
is mapped across multiple vertices of the range tree, in a non-affine fashion.

If $\eta$ is a vertex of $\ut{n}$ we can now define a forward and
backward orbit $O_{\eta}$ of $\eta$ in $\ut{n}$ under the action of
$\langle \alpha\rangle$, to the extent that we restrict ourselves to powers
of $\alpha$ that take the interval under $\eta$ affinely bijectively to
an interval under a vertex of $\ut{n}$.  In more general circumstances,
given any integer $k$ and a vertex $\eta$ of $\ut{n}$, we will use the
notation $\eta \alpha^k$ to denote the subset of $\CS_n$ which is the image
of the interval under $\eta$ under the map $\alpha^k$, and if that set
happens to be the interval under a vertex $\tau$ in $\ut{n}$ (so that
the restriction of $\alpha^k$ to the interval under $\eta$ takes the
interval underlying $\eta$ affinely to the interval underlying $\tau$),
then we may denote the vertex $\tau$ by $\eta \alpha^k$ as well.

\newcommand{\labs}{L_{(A,B,\sigma)}}

Let $L_{(A,B,\sigma)}$ denote the set of vertices of $\ut{n}$ which are
either leaves of $A$ or leaves of $B$, and let $\leafl\in\labs$.

It is possible that $\leafl$ is a neutral leaf whose vertex in $\ut{n}$
has orbit $O_{\leafl}$ entirely contained in the neutral leaves of $A$
and $B$.  By Remark \ref{neutralTravelers}, and the fact that the
neutral leaves are finite in number, we see that in this case the
orbit of $\leafl$ is periodic in $\ut{n}$.  In this case, we call
$\leafl$ a \emph{periodic leaf}.

Now suppose $\leafl$ is not a periodic leaf of $A$.  This implies that
if we consider the forward and backward orbits of $\leafl$ under the
action of $\langle \alpha\rangle$, then in both directions, the orbit will
exit the set of neutral leaves.  In particular, there is a minimal integer
$r\leq 0$ and a maximal integer $s \geq 0$ so that for all integers
$i$ with $r\leq i\leq s$ we have that $\leafl \alpha^i\in L_{(A,B,\sigma)}$.
It is also immediate that $\leafl \alpha^r$ is a leaf of $A-B$ and $\leafl
\alpha^s$ is a leaf of $B-A$, while for all values of $i$ with $r<i<s$ we
see that $\leafl \alpha^i$ is a neutral leaf.

Thus, we have the following lemma.

\bl
\label{IACs}

Suppose $\alpha\in V$ is non-trivial and $\alpha$ is represented by a tree pair
$(A,B,\sigma)$.  If $\leafl\in\labs$, then either $\leafl$ is 
\be
\item a periodic neutral leaf, in which case there is a maximal
  integer $s\geq 0$ so that the \emph{iterated augmentation chain} defined by $IAC(\leafl) := (\leafl \alpha^i)_{i = 0}^s$ is a
  sequence of neutral leaves so that $s+1$ is the smallest positive
  power so that $\leafl \alpha^{s+1} = \leafl$,
\item a leaf of $A-B$, in which case there is a maximal
integer $s>0$ so that $IAC(\leafl) := (\leafl \alpha^i)_{i = 0}^s$ is a
sequence of leaves of $A$ or $B$, and furthermore, we then have
$\leafl \alpha^s$ is a non-neutral leaf of $B$ while $\leafl \alpha^i$ is a
neutral leaf in $\labs$ for all indices $i$ with $0<i<s$,
\item a leaf of $B-A$, in which case there is a minimal
integer $r<0$ so that $IAC(\leafl) := (\leafl \alpha^i)_{i = r}^0$ is a
sequence of leaves of $A$ or $B$, and furthermore, we then have
$\leafl \alpha^r$ is a non-neutral leaf of $A$ while $\leafl \alpha^i$ is a
neutral leaf in $\labs$ for all indices $i$ with $r<i<0$, or
\item a neutral, non-periodic leaf, in which case, $\leafl$ is a
neutral leaf in a sequence $IAC(\eta)$ for some vertex $\eta$ which is a
leaf of $A-B$, as discussed in point $(2)$.  In this case we set $IAC(\leafl):=IAC(\eta)$.
\ee 
\el

We now define and comment on some language from \cite{brinHigherV} and from
\cite{salazarPaper}. Suppose we have the hypotheses of Lemma
\ref{IACs}.
The definition of iterated augmentation chain $IAC(\leafl)$ in Lemma \ref{IACs} 
reflects the fact that we can augment the trees $A$ and
$B$ at each vertex along an iterated augmentation chain, and end up
with a new representative tree pair for $\alpha$ (the first
vertex in an augmentation chain of type $(2)$ or $(3)$ can only be
augmented in the domain tree $A$, while the last such vertex can only
be augmented in the range tree $B$).

If $\leafl$ is of type $(2)$, with $IAC(\leafl) = (\leafl \alpha^i)_{i =
0}^s$ for some positive integer $s$, where $\leafl \alpha^s$ is an ancestor of
$\leafl$, then we say that \emph{$\leafl$ is a repeller}.  If $\leafl$
is of type (3), with $IAC(\leafl) = (\leafl \alpha^i)_{i = r}^0$ for some
negative integer $r$, and where $\leafl \alpha^r$ is an ancestor of $\leafl$, then
we say that \emph{$\leafl$ is an attractor}.

We are now ready to define what it means for a tree pair to be a
revealing pair.  Suppose $\alpha\in V$ and the tree pair $(A,B,\sigma)$
represents $\alpha$.  If every component of $A-B$ contains a repeller, and
every component of $B-A$ contains an attractor, then we say that
$(A,B,\sigma)$ is a \emph{revealing pair representing $\alpha$}.

The discussion beginning section 10.7 in \cite{brinHigherV} proves
that every element of $V_n$ admits a revealing pair.  It is not hard to
generate an algorithm which will transform any representative tree
pair for an element of $V_n$ into a revealing pair.

Given $\alpha\in V_n$, we will denote by $R_\alpha$ the set of all revealing pairs
for $\alpha$. We will use the symbol $\sim$ as a relation in the
fashion $(A,B,\sigma)\sim \alpha$ denoting the fact that $(A,B,\sigma)\in
R_\alpha$.  In this case, we can further name leaves of $A-B$ and $B-A$.
If $\leafl$ is a leaf of $A-B$ and $\leafl$ is not a repeller, then we
say $\leafl$ is a \emph{source}.  If $\leafl$ is a leaf of $B-A$ and $\leafl$
is not an attractor, then we say $\leafl$ is a \emph{sink}.

There are finitely many process types called ``Rollings'' introduced
by Salazar in \cite{salazarPaper}.  Rollings are methods by which one
can carry out a finite collection of simple expansions to a revealing
tree pair $(A,B,\sigma)$ to produce a new revealing pair
$(A',B',\sigma')$.   

Below, we give the definitions and some discussion for rollings of type II.  We give definitions for the other types of rollings in sub-section \ref{discreteTrainTracks}.

The tree pair $(A',B',\sigma')$ is a \emph{single rolling of type II
  from $(A,B,\sigma)$} if it is obtained from $(A, B,\sigma)$ by adding
a copy of a component $U$ of $A-B$ to $A$ at the last leaf in the
orbit of the repeller in $U$ and to $B$ at its image; or, by adding a
copy of a component $W$ of $B-A$ to $A$ at the first leaf in the
orbit of the attractor (the leaf of $A$ corresponding to the root node of $W$ in $B$) and to $B$ at its image.

The tree pair $(A',B',\sigma')$ is a \emph{rolling of type II from
  $(A,B,\sigma)$} if it is obtained from $(A,B,\sigma)$ by a finite
collection of single rollings of type II applied to the initial tree
pair $(A,B,\sigma)$ in some order.

We now state three lemmas about properties of revealing pairs.  All of
these properties are fairly straightforward to verify.  In the cases
of Lemma \ref{RepellerProps} and Lemma \ref{AttractorProps}, the
curious reader may also refer to the discussion in Sections 3.3 and 3.4
of \cite{salazarPaper} for alternative proofs.

Below, we slightly abuse the notion of a vertex name, by
associating an infinite descending path with an infinite ``name''
string. This then represents a point in the Cantor set which is the boundary of the infinite tree.

Our first lemma discusses repellers, sources, and fixed points. 

\bl 
\label{RepellerProps} 

Suppose $\alpha\in V_n$ and $(A,B,\sigma)\sim \alpha$, and that  $\leafl$ is a leaf of
a component $C$ of $A-B$, so that $IAC(\leafl) = (\leafl \alpha^i)_{i = 0}^s$.

\be
\item If $\leafl$ is a repeller (so $\leafl \alpha^s$ is an ancestor of
  $\leafl$), set $\gamma_i$ to be the name of the node $\leafl \alpha^i$
  for each $0\leq i\leq s$. Set $\Gamma$ to be the string which is the suffix
  one would append to the name of the node $\gamma_s$ to obtain the name of the node $\gamma_0$
  i.e. the path from the ancestor to the repeller. We notate this by the expression
  $\gamma_s\Gamma =\gamma_0$. Thus the infinite
  descending path corresponding to the infinite string
  $\gamma_i(\Gamma)^{\infty}$ represents a unique repelling fixed
  point $p_{\gamma_i}$ in the interval $X_i$ of $\CS_n$ underneath
  $\gamma_i$, under the action of $\langle \alpha^s\rangle$, for each index
  $i$ with $0\leq i<s$.

\item if $\leafl$ is a source, then $\leafl \alpha^s$ is a sink.
\item $\langle \alpha\rangle \cong \Z$.  
\ee 
\el 

{\it Proof:} 

The inverse $\alpha^{-s}$ of $\alpha^s$ maps the Cantor set $X_s$ under
$\gamma_s$ bijectively to the Cantor set $X_0$ under $\gamma_0$, where
$X_0\subset X_s$, in affine fashion.  The only infinite descending
path which is fixed by this map is the path terminating in
$\Gamma^{\infty}$, thus $p_{\gamma_0} = \gamma_0\Gamma^{\infty} =
\gamma_s\Gamma\Gamma^{\infty}$ is the unique attracting fixed point of
$\alpha^{-s}$ in $X_s$, and is thus the unique repelling fixed point of
$\alpha^s$ within $X_0$ (and even within all of $X_s$) under the action of
$\alpha^s$.

We now show below why each of the points
$p_{\gamma_i}$ are also repelling fixed points of $\alpha^s$.

We first obtain a new revealing tree pair $(A',B',\sigma')\sim \alpha$ which
is a single rolling of type II from $(A,B,\sigma)$, constructed as
follows. 

  Glue a copy of $C$ at $\leafl_s$ in $B$ to produce $B'$, and a
  further copy of $C$ at the leaf $\gamma_{s-1}$ of $A$ to produce
  $A'$.  The resulting tree pair $(A',B',\sigma')$ so obtained has all
  of the leaves of $C$ of $A-B$ as neutral leaves (except in the case
  where $s = 1$, in which case $A'$ is $A$ with a copy of $C$ attached
  at $\gamma_0$, a leaf of the original $C$).

In any case, the original copy of $C$, within $A'$, is now contained
in $A'\cap B'$.  The new copy of $C$ in $A'$ is a complementary
component of $A'-B'$ and contains the leaf $\gamma_{s-1}\Gamma$ as a
repeller.  There are no other new complementary components for the
tree pair $(A', B',\sigma')$, which therefore must represent a
revealing pair for $\alpha$.

Now, by a minor adjustment to the argument in first paragraph, the
point $p_{\gamma_{s-1}}$ is a fixed repelling point of $\alpha^s$.

We can now continue inductively in this fashion to show that each
point $p_{\gamma_i}$ is the unique repelling fixed point of $X_i$
under the action of $\alpha^s$ by building a revealing pair
$(A'',B'',\sigma'')$ for $\alpha$ with the point $p_{\gamma_i}$ as a point
in a repelling leaf of a complementary component of
$(A'',B'',\sigma'')$ with shape $C$ and spine $\Gamma$ rooted at
$\gamma_i$.

We leave the second point of the lemma to the reader, while the third
point is immediate from the first since some power of $\alpha$ has a
repelling fixed point.

\qquad$\diamond$

We call each subset $X_i$ a \emph{basin of repulsion for $\alpha$}, since all points in $X_i$
eventually flow out of $X_i$ under repeated iteration of $\alpha^s$, never to return, except
for $p_{\gamma_i}$, for each $0\leq i\leq s$.  We call the string $\Gamma$
above the \emph{spine of the repeller $\gamma_0$} or the \emph{spine
  of $C$}, as it describes the shape of the path in $C$ from the root
of $C$ to the repeller $\gamma_0$.

We call each point $p_{\gamma_i}$ a \emph{periodic repelling point of
  $\alpha$} for each index $i$ with $0\leq i<s$.  We may also call these
points \emph{fixed repelling points of $\alpha^s$}, noting in passing that
not all periodic repelling points of $\alpha$ are necessarily fixed by
$\alpha^s$.  In similar fashion we call the sequence of points
$(p_{\gamma_i})_{i = 0}^{s-1}$ a \emph{periodic orbit of periodic
  repelling points for $\alpha$}.  We denote by $\mathcal{R}_\alpha$ the set of
periodic repelling points of $\alpha$, noting that it is a finite set.

We now give a similar lemma discussing attractors and sinks.

\bl \label{AttractorProps} 

Suppose $\alpha\in V_n$ and $(A,B,\sigma)\sim \alpha$ and $\leafl$ is a leaf of
a component $C$ of $B-A$, so that $IAC(\leafl) = (\leafl \alpha^i)_{i =
  r}^0$, for some negative integer $r$.

\be

\item If $\leafl$ is an attractor (so $\leafl \alpha^r$ is an ancestor of
  $\leafl$), set $\gamma_i$ to be the name of the node $\leafl \alpha^i$
  for each $r\leq i\leq 0$.  In this case the string $\Gamma$ which
  has $\gamma_r\Gamma =\gamma_0$ has the property that the infinite
  descending path corresponding to the name
  $\gamma_i(\Gamma)^{\infty}$ represents a unique attracting fixed
  point $p_{\gamma_i}$ in the interval $X_i$ of $\CS_n$ underneath
  $\gamma_i$, under the action of $\langle \alpha^r\rangle$, for each index
  $i$ with $r\leq i <0$.

\item If $\leafl$ is a sink, then $\leafl \alpha^r$ is a source.
\item $\langle \alpha\rangle \cong \Z$.
 
\ee

\el 

{\it Proof:} 

This proof is similar to the proof of the previous lemma, where here
$\alpha^{-1}$ has a tree pair with $\leafl$ as a repeller, and each
$p_{\gamma_i}$ is a periodic repelling point of $\alpha^{-1}$.

\qquad$\diamond$

We now extend the notation from Lemma \ref{RepellerProps} to apply to the sets named in Lemma \ref{AttractorProps} as below.

We call each subset $X_i$ a \emph{basin of attraction for $\alpha$
  as indicated by $(A,B,\sigma)$}, since all points in $X_i$
eventually limit to $p_{\gamma_i}$ under repeated iteration of $\alpha^{-r}$, for each $r\leq i\leq 0$.  We call the string $\Gamma$
above the \emph{spine of the attractor $\gamma_0$} or the \emph{spine
  of $C$}, as it describes the shape of the path in $C$ from the root
of $C$ to the attractor $\gamma_0$.

We call each point $p_{\gamma_i}$ a \emph{periodic attracting point of
  $\alpha$} for each index $i$ with $r\leq i\leq 0$.  We may also call
these points \emph{fixed attracting points of $\alpha^{-r}$}.  In similar
fashion we call the sequence of points $(p_{\gamma_i})_{i = r}^{-1}$
a \emph{periodic orbit of periodic attracting points for $\alpha$}.  We
denote by $\mathcal{A}_\alpha$ the set of periodic attracting points of $\alpha$,
noting that it is a finite set.

In the previous two lemmas, if $\leafl$ is a source or a sink (case
two in each lemma), we refer to $IAC(\leafl)$ as a \emph{source-sink
  chain}.

The next lemma follows directly from the two above, and the
classification of the orbits of the leaves of $A$ and $B$.  This lemma is 
a version of one result proved by Burillo, Cleary, Stein and Taback in their joint
work \cite{BCSTCombinatoricsT} as Proposition 6.1.  We give
a new proof here, as the situation is greatly simplified through the
use of revealing pairs.

\bl \label{torsionRevPair}

Suppose $\alpha\in V_n$.  There is an integer $n$ so that $\alpha$ has order $n$
if and only if there is a tree pair $(A,B,\sigma)$ representing $\alpha$
with $A = B$.

\el

{\it Proof:}

Suppose $\alpha$ does not have infinite order, and let $(A,B,\sigma)\sim \alpha$
be a revealing pair representing $\alpha$.  We must have that
$(A,B,\sigma)$ admits no repellers or attractors, thus both $A-B$ and
$B-A$ are empty, and so $A = B$.

Suppose instead $(A,B,\sigma)$ is a tree pair representing $\alpha$ with $A
= B$. Then it is immediate from the definition of multiplication for
tree pairs that the order of $\alpha$ is the order of the permutation
$\sigma$.  \qquad $\diamond$

We denote by $\mathcal{P}_\alpha$ the points of $\CS_n$ which underlie the periodic neutral leaves of $\labs$. 
Note that the set $\mathcal{P}_\alpha$ is independent of the choice of revealing pair used to represent $\alpha$. We further denote by $Per(\alpha)$ the set of all periodic points of $\alpha$, that is 
\[
Per(\alpha) = \mathcal{R}_\alpha\sqcup\mathcal{A}_\alpha\sqcup\mathcal{P}_\alpha.
\]

The following tree pairs represent the previously defined element $\theta\in T$ (we apply
an augmentation to our first tree pair to produce a revealing pair
representing $\theta$).
\begin{center}
\psfrag{1}[c]{$1$}
\psfrag{2}[c]{$2$}
\psfrag{3}[c]{$3$}
\psfrag{4}[c]{$4$}
\psfrag{5}[c]{$5$}
\psfrag{6}[c]{$6$}
\psfrag{7}[c]{$7$}
\psfrag{D}[c]{$D$}
\psfrag{R}[c]{$R$}
\includegraphics[height=165pt,width=240 pt]{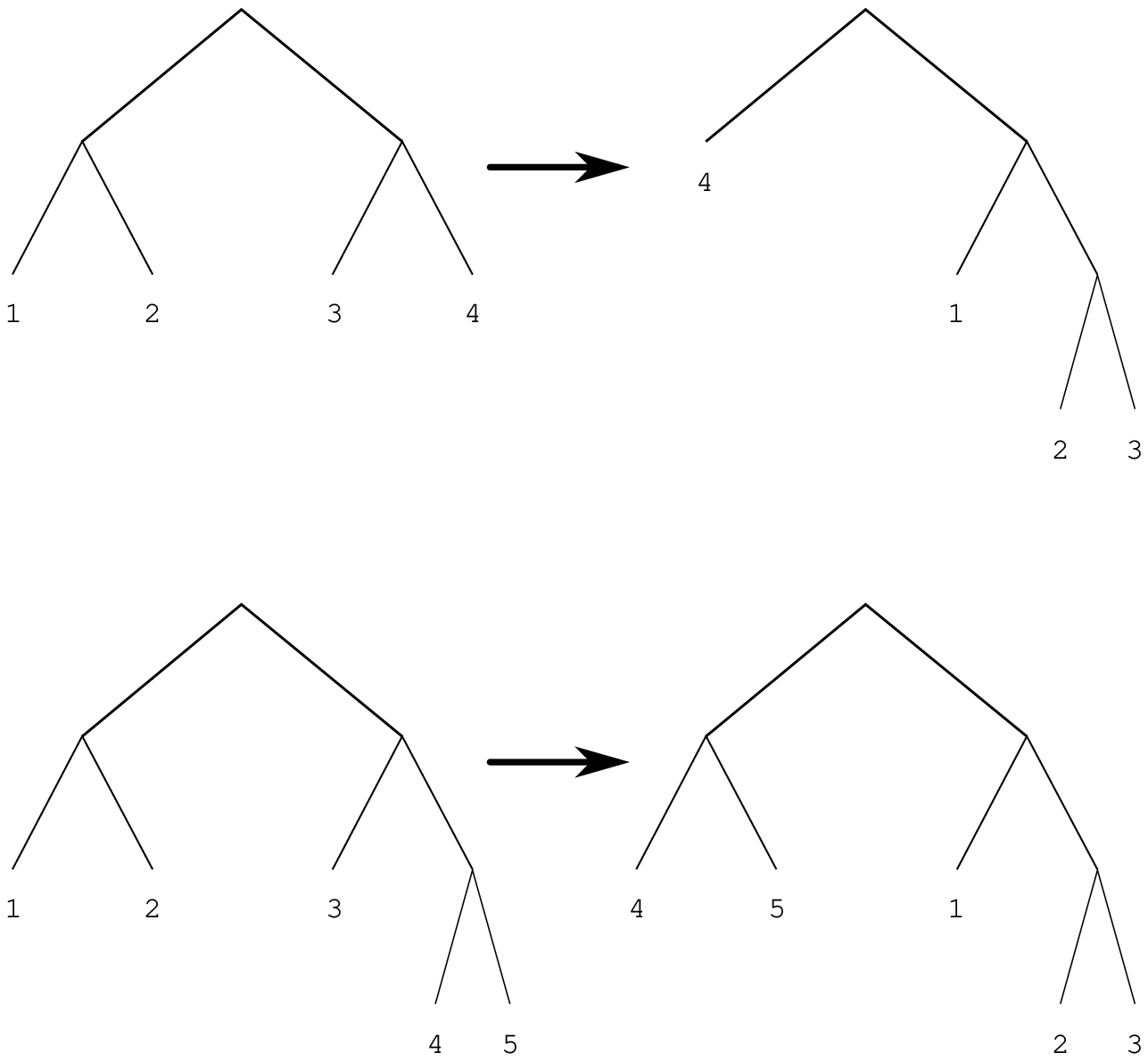}
\end{center}
A quick examination of the second tree pair above should convince the
reader that $\theta$ has order five and rotation number $2/5$.

 We now give a series of diagrams for a revealing pair $(A,B,\sigma)$
 representing a particular non-torsion element $\alpha\in V$.  This element
 has a revealing pair which contains many of the structures we have
 been discussing.  In each diagram below we illustrate some of the
 particular aspects we have discussed above.

Below is an example of our revealing pair $(A,B,\sigma)$, with the
neutral leaves underlined.  Recall that the left tree is $A$, and the right tree is $B$.
\begin{center}
\psfrag{1}[c]{$\,\,1$}
\psfrag{2}[c]{$\,\,2$}
\psfrag{3}[c]{$\,\,3$}
\psfrag{4}[c]{$\,\,4$}
\psfrag{5}[c]{$\,\,5$}
\psfrag{6}[c]{$\,\,6$}
\psfrag{7}[c]{$\,\,7$}
\psfrag{8}[c]{$\,\,8$}
\psfrag{9}[c]{$\,\,9$}
\psfrag{10}[c]{$\,\,10$}
\psfrag{D}[c]{$D$}
\psfrag{R}[c]{$R$}
\includegraphics[height=150 pt,width=300 pt]{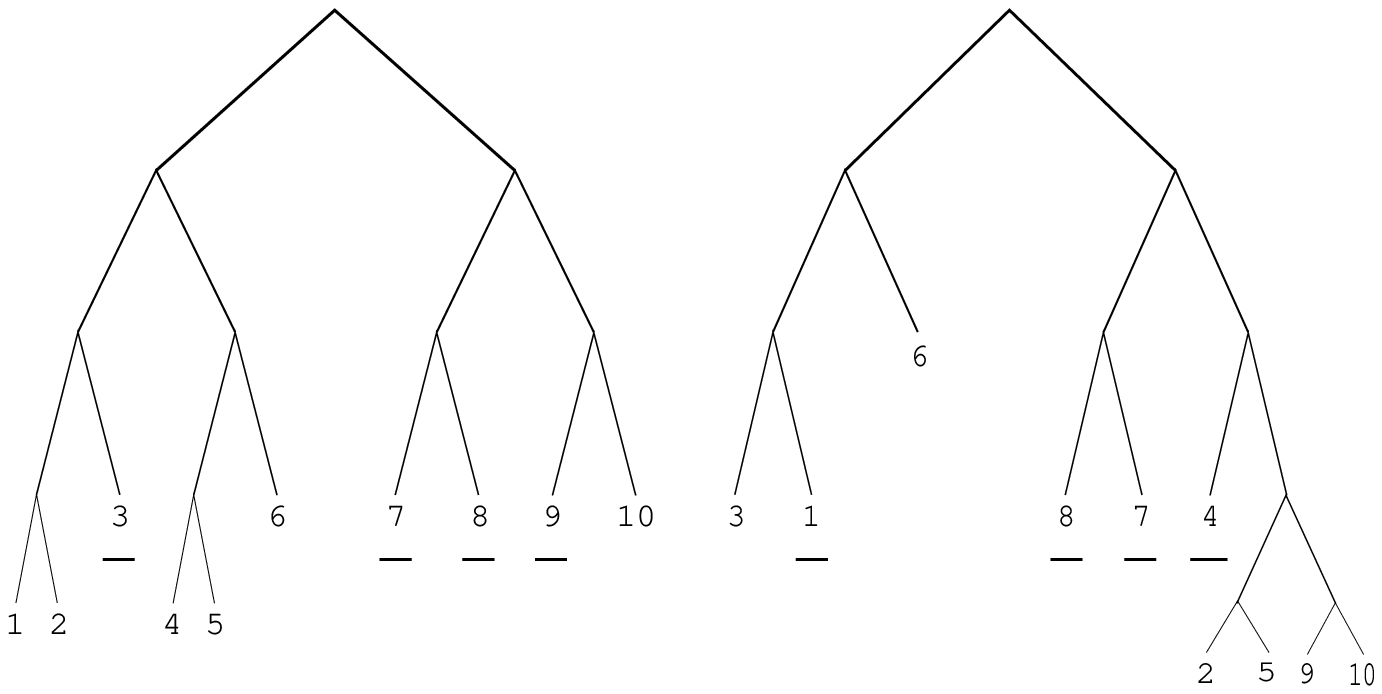}
\end{center}
In the next diagram, we point out a periodic orbit of neutral leaves of length two.
\begin{center}
\psfrag{1}[c]{$\,\,1$}
\psfrag{2}[c]{$\,\,2$}
\psfrag{3}[c]{$\,\,3$}
\psfrag{4}[c]{$\,\,4$}
\psfrag{5}[c]{$\,\,5$}
\psfrag{6}[c]{$\,\,6$}
\psfrag{7}[c]{$\,\,7$}
\psfrag{8}[c]{$\,\,8$}
\psfrag{9}[c]{$\,\,9$}
\psfrag{10}[c]{$\,\,10$}
\psfrag{D}[c]{$D$}
\psfrag{R}[c]{$R$}
\includegraphics[height=150 pt,width=300 pt]{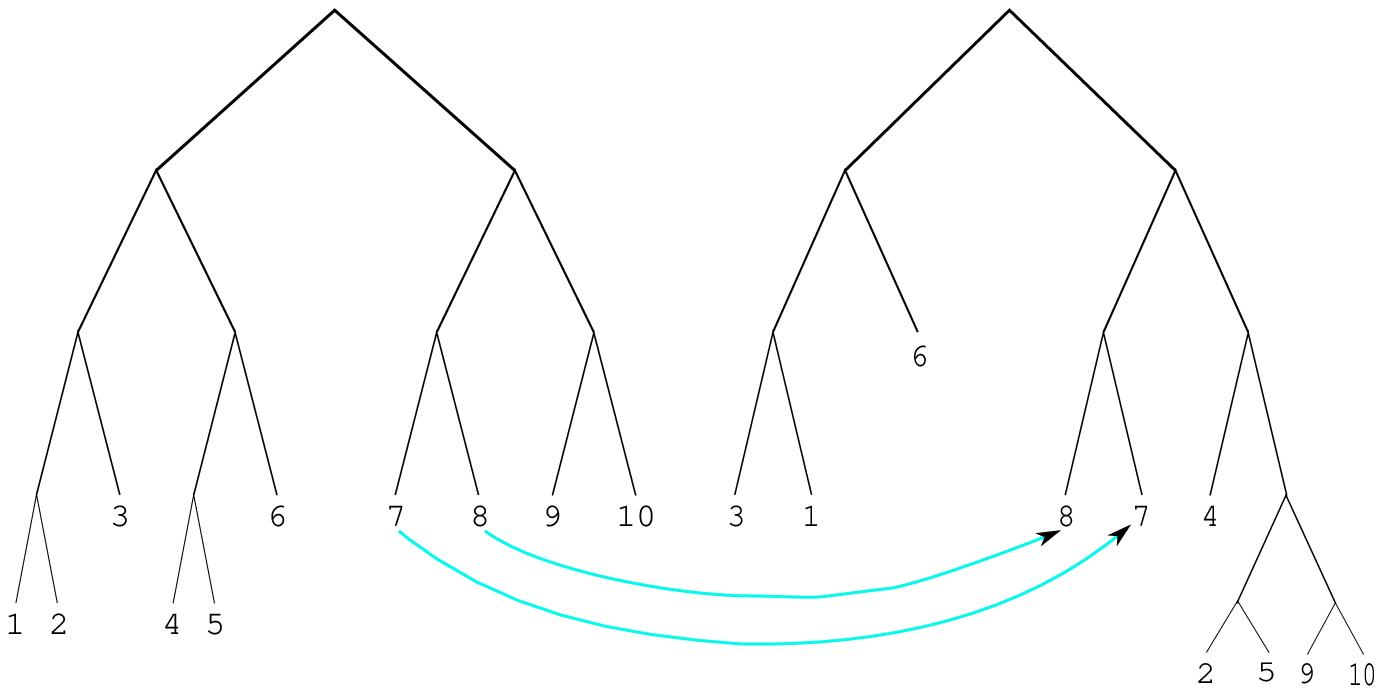}
\end{center}
Below, the vertex $11111$ in $B$ is an attractor; its
iterated augmentation chain is $(111,11111)$, a sequence of length
two.  The fixed point of the attractor corresponds to the value $1$ in
the unit interval.  The word $\Gamma$ for this attractor is $11$.
\begin{center}
\psfrag{1}[c]{$\,\,1$}
\psfrag{2}[c]{$\,\,2$}
\psfrag{3}[c]{$\,\,3$}
\psfrag{4}[c]{$\,\,4$}
\psfrag{5}[c]{$\,\,5$}
\psfrag{6}[c]{$\,\,6$}
\psfrag{7}[c]{$\,\,7$}
\psfrag{8}[c]{$\,\,8$}
\psfrag{9}[c]{$\,\,9$}
\psfrag{10}[c]{$\,\,10$}
\psfrag{D}[c]{$D$}
\psfrag{R}[c]{$R$}
\includegraphics[height=150 pt,width=300 pt]{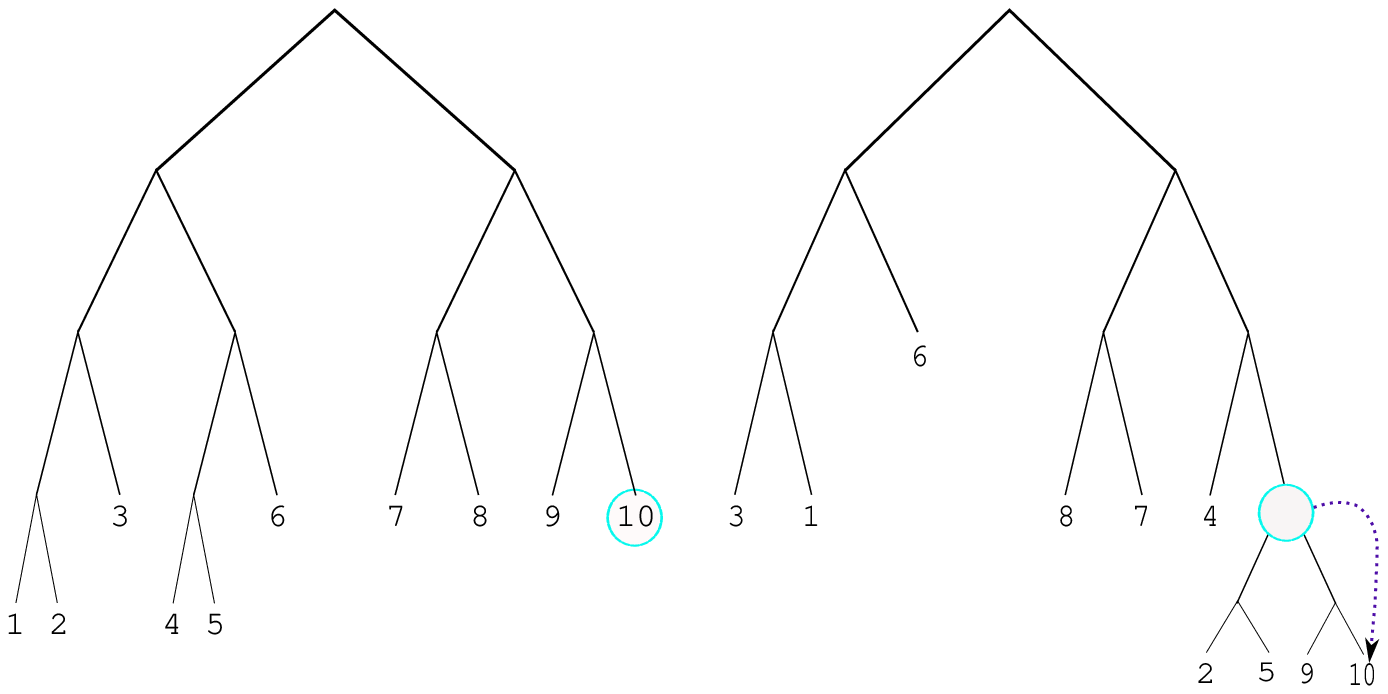}
\end{center}

In the diagram to follow, the vertices $0000$ and $011$ represent
repellers in $A$.  The dotted paths track the orbits of the
repellers along their iterated augmentation chains.
\begin{center}
\psfrag{1}[c]{$\,\,1$}
\psfrag{2}[c]{$\,\,2$}
\psfrag{3}[c]{$\,\,3$}
\psfrag{4}[c]{$\,\,4$}
\psfrag{5}[c]{$\,\,5$}
\psfrag{6}[c]{$\,\,6$}
\psfrag{7}[c]{$\,\,7$}
\psfrag{8}[c]{$\,\,8$}
\psfrag{9}[c]{$\,\,9$}
\psfrag{10}[c]{$\,\,10$}
\psfrag{D}[c]{$D$}
\psfrag{R}[c]{$R$}
\includegraphics[height=150 pt,width=300 pt]{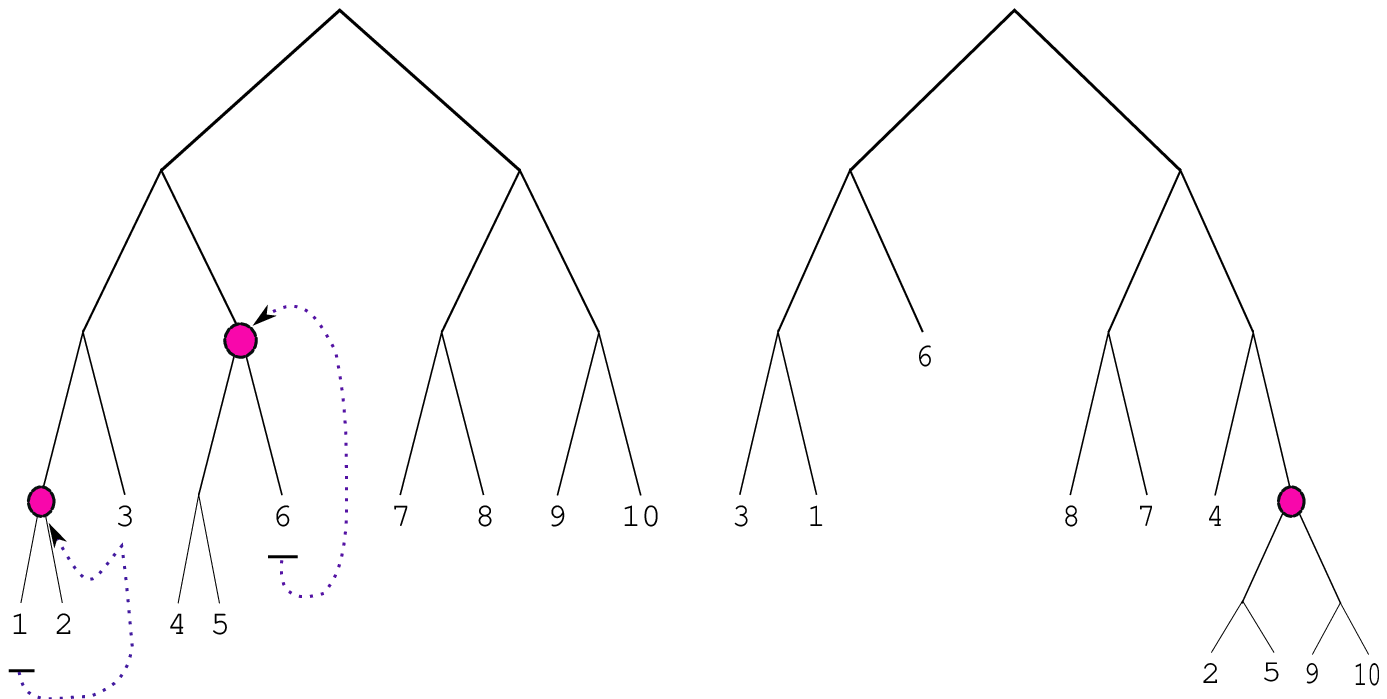}
\end{center}
Finally, we highlight the fact that sources flow to sinks.
Note how the lengths of the paths from sources to sinks are not
uniform.  In particular, the source $0100$ first hops to $110$
before next landing in the basin of attraction under the vertex $000$.
\begin{center}
\psfrag{1}[c]{$\,\,1$}
\psfrag{2}[c]{$\,\,2$}
\psfrag{3}[c]{$\,\,3$}
\psfrag{4}[c]{$\,\,4$}
\psfrag{5}[c]{$\,\,5$}
\psfrag{6}[c]{$\,\,6$}
\psfrag{7}[c]{$\,\,7$}
\psfrag{8}[c]{$\,\,8$}
\psfrag{9}[c]{$\,\,9$}
\psfrag{10}[c]{$\,\,10$}
\psfrag{D}[c]{$D$}
\psfrag{R}[c]{$R$}
\includegraphics[height=150 pt,width=300 pt]{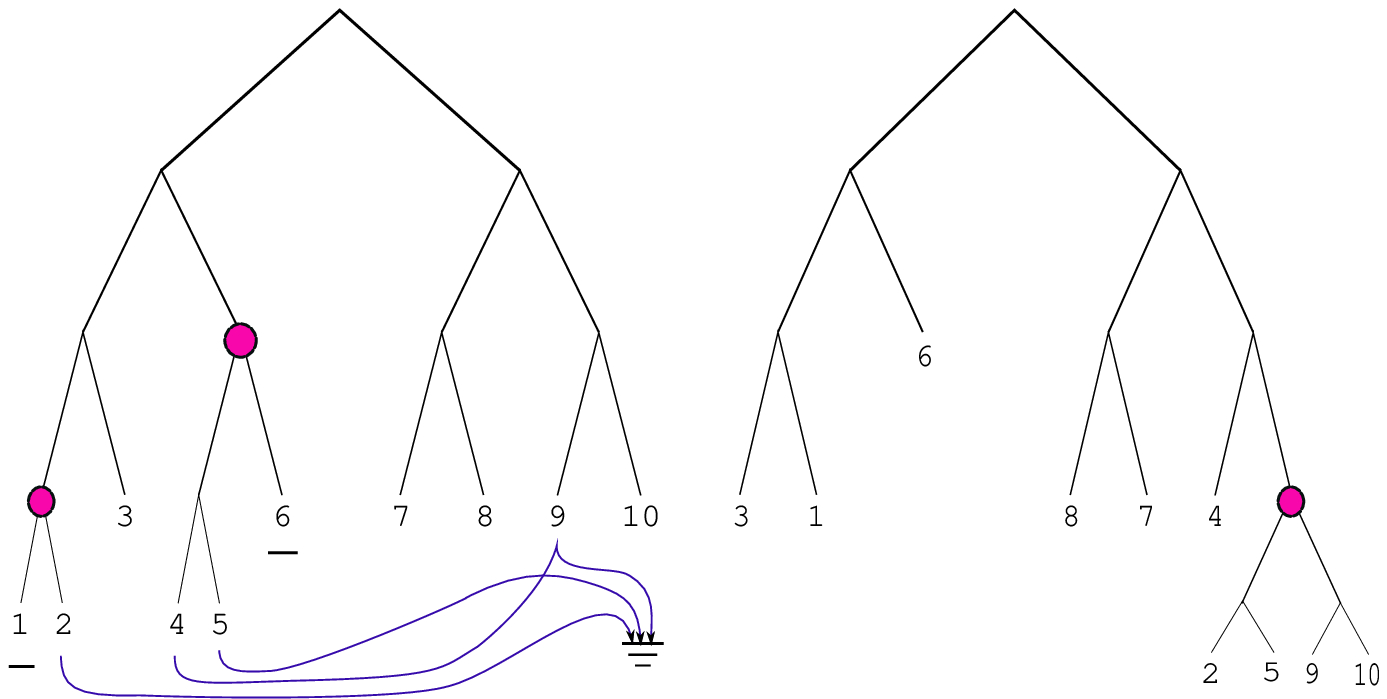}
\end{center}

\subsection{Discrete train tracks and laminations \label{discreteTrainTracks}}
In this subsection we are concerned with modelling the dynamics in $\mathfrak{C}_n$ under the action of the group generated by an element $v$ of $V_n$.  Naturally, the relevant information is contained in any particular revealing pair $\rho = (C,D,\theta)$ for the element, but we have found that tools from surface topology afford us another way to visualize these dynamics.  In particular, these dynamics can be usefully described by a combinatorial object introduced by Thurston  (see \cite{Thurston1, Thurston2}) to aid in the study of surfaces, namely, a train track.  Given a revealing pair, it is easy to draw a combinatorial train track which in turn ``carries" a lamination in some compact surface with boundary.  It is this lamination which models, in some sense, the movement of points in the Cantor set under iteration of the map $v$.  As in the theory of laminations carried by a train track, one quickly realizes that most of the relevant visual information is actually contained in the train-track object (see for instance, Penner and Harer's book \cite{PennerHarer}).  We call the train track object developed here a \emph{discrete train track}, even though it is continuous in nature; the name is meant to emphasize that our train tracks model dynamics in a totally disconnected set under iterations of a fixed map.  

  Below, we describe in detail our method for generating a discrete train track $TT_\rho$ from the revealing pair $\rho$ representing the general element $v \in V_n$ mentioned above, and briefly describe how to model the lamination it carries.   We build an example discrete train track and lamination using the element $\alpha$ and the tree pair $(A,B,\sigma)$ from the previous subsection.  Finally, we discuss some of the utility of $TT_\rho$ and some facts about how $TT_\rho$ would change under basic operations applied to $v$ (conjugation, then a choice of representative revealing pair) or $\rho$ (Salazar's rollings).    In the next subsection, we describe a derived object, the flow graph, which carries less information than the train track (although flow graphs in general still contain enough information to answer many dynamics questions).  In our current practice, we find discrete train tracks and flow graphs to be helpful for understanding dynamics, while the generating revealing pairs tend to be helpful for any involved computations and for specifying elements with appropriate desired dynamics.  

Note that any such discrete train track $TT_\rho$ is a representative of an equivalence class of similar train tracks, where the equivalence class is determined by all the discrete train tracks for revealing pairs equivalent to $\rho$ up to Salazar's rollings, choice of location for drawing the complementary trees along the orbits of repelling and attracting fixed points, and conjugacy in $V_n$.  Thus we are picking a representative object which is in some sense less well chosen than some ``minimal" train track in this class.  The result of applying Belk and Matucci's geometric conjugacy invariant, the closed strand diagram, is very close to what a ``minimal" discrete train track for an element of $V_n$ should be (leaving out the demarcations representing iterating the element, and of course, including some simplifications which result from the effects of conjugation (see \cite{BelkMatucci})).

Here is how one draws a discrete train track from a revealing pair.

\begin{enumerate}
\item List all of the iterated augmentation chains for the revealing pair.
\item For each chain representing a non-trivial orbit of a periodic neutral leaf:
\begin{enumerate} 
\item Draw a circle.
\item If the chain represents an orbit of length $r$, demarcate the circle into $r$ equal subintervals (typically we end these intervals with dots), and label these with the names of the nodes carrying the neutral leaves, in a clockwise order.
\end{enumerate}
\item For each chain of a repeller:
\begin{enumerate} 
\item Draw a circle.
\item If the repelling periodic point travels an orbit of length $r$, then demarcate the circle into $r$ equal subintervals as above, and label these with the names of the nodes carrying the orbit of the repeller in a counter-clockwise order.
\item For the segment corresponding to the interval of the repelling periodic point of the complementary tree, instead of the one label mentioned in the last point, we label the two ends of the segment with the top and bottom of the spine of the repeller (top (root node of complementary component) before bottom in counter-clockwise order).
\item Lay the complementary component of the repeller along the sub-arc of the circle with the spine labels, gluing the spine to the circle. 
\begin {enumerate}
\item Scale the complementary component so that the spine is the length of the appropriate sub-arc of the circle.
\item Smooth out the tree (and the spine in particular) and bend the spine so that the spine has the same shape as the sub-arc of the circle with the labels from the spine nodes (preserving the current length of the spine).
\item Rigidly rotate the scaled, smoothed, and bent tree in the plane, and translate it so that the spine can be identified with the appropriate sub-arc of the circle. 
\end{enumerate}
\end{enumerate}
\item For each chain for an attractor:
\begin{enumerate} 
\item Draw a circle.
\item If the attracting periodic point travels an orbit of length $r$, then demarcate the circle into $r$ equal subintervals, and label these with the names of the nodes carrying the orbit of the attractor in a counter-clockwise order.
\item For the segment corresponding to the interval of the attracting periodic point of the complementary tree, instead of the one label mentioned in the last point, we label the two ends of the segment with the top and bottom of the spine of the attractor (bottom (attracting leaf node of complementary component) before top in counter-clockwise order). 
\item Lay the complementary component of the attractor along the sub-arc of the circle with the spine labels, gluing the spine to the circle. 
\begin {enumerate}
\item reflect the complementary component of the attractor across a vertical axis.  Thus, the root of the component is now drawn at the bottom, while the left and right hand sides are preserved, respectively, as left and right hand sides.
\item Scale the complementary component so that the spine is the length of the appropriate sub-arc of the circle.
\item Smooth out the tree (and the spine in particular) and bend the spine so that the spine has the same shape as the sub-arc of the circle with the labels from the spine nodes (preserving the current length of the spine).
\item Rigidly rotate the reflected, scaled, smoothed, and bent tree in the plane, and translate it so that the spine can be identified with the appropriate sub-arc of the circle.
\end{enumerate}
\end{enumerate}
\item For each source-sink chain, drawn a line connecting the appropriate sources and sinks (lines may have to cross each other in the case of $V_n$ which is why the lamination carried by the train track can only be embedded in an surface with boundary; strips can pass ``under'' each other).  Demarcate each line with dots representing the length of the source-sink chain, and label sub-arcs with appropriate node labels as above for other sorts of iterated augmentation chains.
\item Add parenthetical labels for splittings of trees where the support of the whole tree will be mapped away by one application of the element.  Add parenthetical labels anywhere else as desired to improve clarity.  (This step is not strictly necessary, but we find it to be helpful.)
\end{enumerate}

If we follow the process above for the element $\alpha$, the diagram below is an example of what we may obtain (we include some drawn under-crossings following the methods from drawing knots from knot theory).  

\begin{center}
\includegraphics[height=250 pt, width= 450 pt]{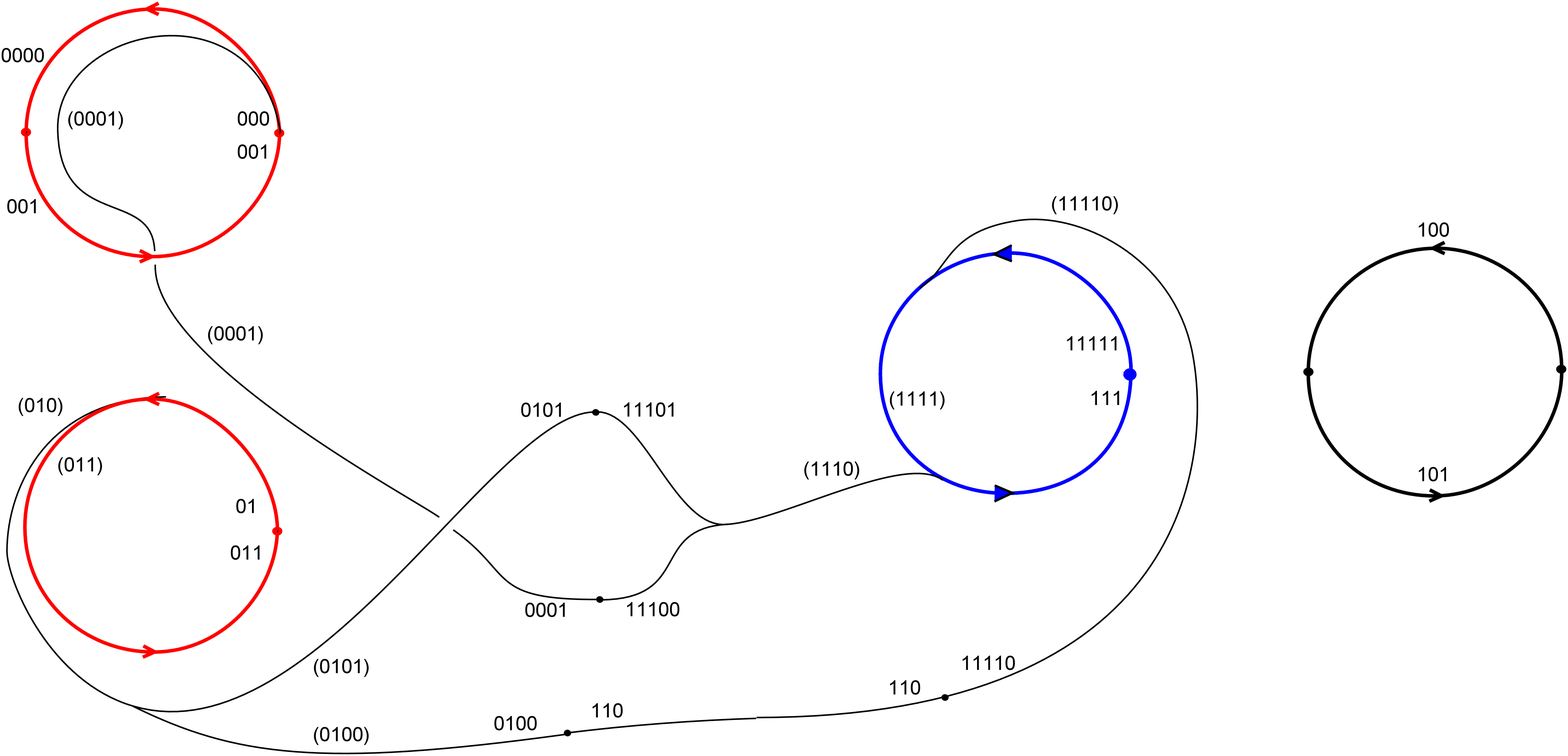}
\end{center}

To visualize the lamination carried by a discrete train track, for each sub-arc in the demarcations of the discrete train tracks, one can draw a transverse copy of an $I$-fiber with an embedded copy of the Cantor set $\CS_n$.  If we stabilize this by building a product with the interval $I$ (which lines we draw parallel to any train track sub-arc), then one has a picture of the lamination carried by the train track. (Glue the ends of the $I$-fibers of the Cantor sets using the rules provided by the map $v$, as a local picture of a mapping cone on the $I$-fiber running transverse to any demarcation dot, allowing contraction and expansion in the transverse $I$-fibers near to any tree-splitting.  

Note that by reflecting the attracting complementary components as directed, we are able to draw the resulting carried lamination on a compact orientable surface with boundary. (The ``left" and ``right" portions of a Cantor set underlying a node are correctly associated without any twisting.)
  
Separately, one can observe that the group homomorphism $\mathcal{S}$ of the penultimate subsection of this article is connected with holonomy measurements for the carried lamination along repelling cycles when applied to the element we are centralizing. 

A local diagram of the lamination carried by this train track local to the repelling cycle in the upper left corner
is included in the diagram below. One can see the rescaling near the periodic repeller and portions of the Cantor
set moving away along flow lines.

\begin{center}
\includegraphics[height=125 pt, width= 400 pt]{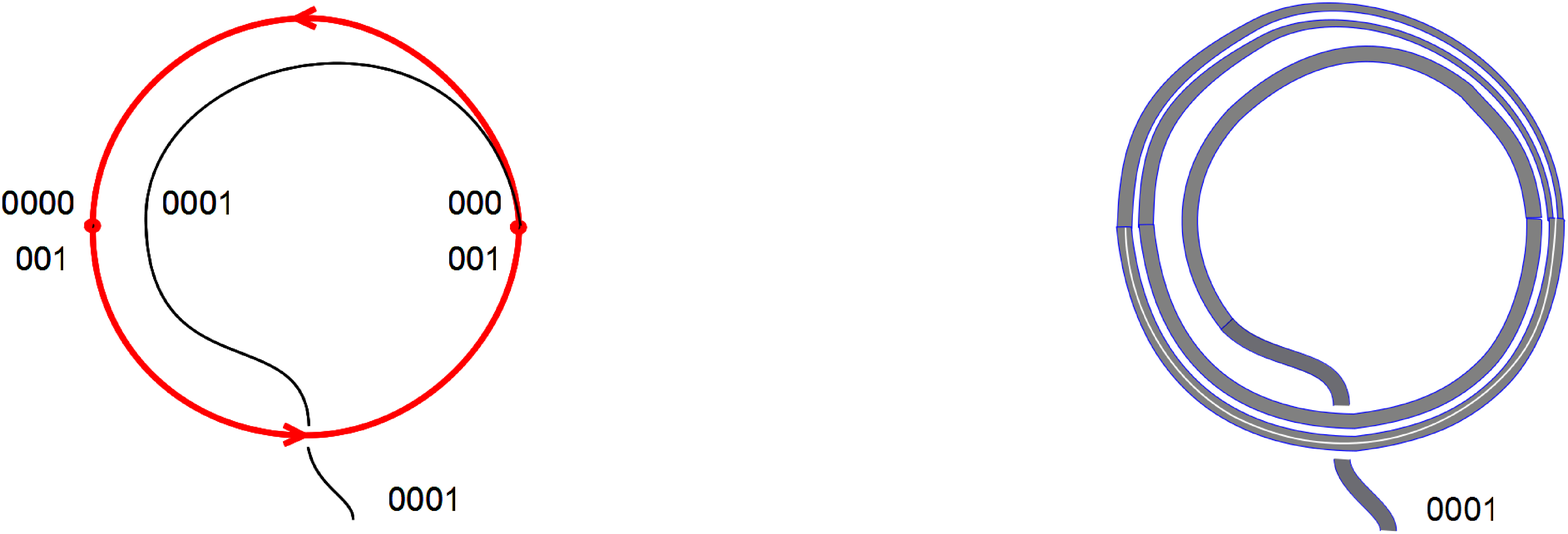}
\end{center}

Let $TT_1$ be a discrete train track drawn in accordance with the method above.  We define the \emph{support of the discrete train track $TT_1$} to be given as the union of the Cantor sets underlying the node labels in the iterated augmentation chains of the (non-trivial) periodic orbits, the source-sink chains, and the orbits of the repellers and attractors used to create the drawing of $TT_1$.  Seeing a discrete train track as a (possibly disconnected) graph, we define the components of $TT_1$ to be the connected components of $TT_1$ as a graph.  Given any such component $X$, we describe the \emph{support of the component $X$} to be the union of the Cantor sets underlying the node labels of $X$.  We call a component a \emph{torsion component} if the component is a circle generated from a periodic orbit of neutral leaves, otherwise we call a component a \emph{non-torsion component.}  Similar language will be defined for flow graphs, below.

We now mention some technical facts to do with the relationships amongst discrete train tracks drawn from distinct revealing pair representatives for an element of $V_n$, and also relationships which may arise as a consequence of conjugacy.

 We give the basic definitions of Salazar's rollings below, with the intention of understanding of the variance of discrete train tracks across the set of all representative revealing pairs for an element of $V_n$.  The reader is encouraged to read Salazar's Section 3.5 in \cite{salazarPaper}, where she defines rollings of various types and traces the various impacts of rollings on revealing tree pairs.  

The tree pair $(A',B',\sigma')$ is a \emph{single rolling of type E
  from $(A,B,\sigma)$} if it is obtained from $(A, B,\sigma)$ by adding
an $n$-caret to $A$ and $B$ along each of the leaves of an iterated augmentation chain corresponding to one of two types of iterated augmentation chains; either to all of the leaves in $A$ and $B$ of a periodic orbit of neutral leaves, or to the initial source leaf in $A$ for a source-sink chain, and then to all of the neutral leaves in $A$ and $B$ of that chain, and then to the leaf of $B$ corresponding to the sink of that chain.  (These are called \emph{elementary rollings}.)

The tree pair $(A',B',\sigma')$ is a \emph{single rolling of type I
  from $(A,B,\sigma)$} if it is obtained from $(A, B,\sigma)$ by adding a cancelling tree along all of the leaves of $A$ along the orbit of a repeller, and at $B$ at the image of these leaves under the map, or by adding a cancelling tree at all the leaves of $A$ which appear in the reverse orbit of an attractor, and at the leaves of $B$ to which these leaves of $A$ are mapped.   If $W$ is the complementary component of $A-B$ or $B-A$ corresponding to the repeller or attractor in this discussion, and $W$ is rooted at node $\Sigma$, then a tree $C$ is a \emph{cancelling tree for $W$} if it is obtained from $W$ by first choosing a proper, non-empty prefix $\Delta$ of the spine $\Gamma$ of $W$ (so that $\Gamma = \Delta\Theta$ for some suffix $\Theta$), and then taking $C$ to be the maximal sub-tree of $W$ which has root $\Sigma$ and containing the node $\Sigma\Delta$ as a leaf.  (Note that if one carries out this process, the corresponding complementary component $W'$ created in $A'-B'$ or $B'-A'$ for the tree pair $(A',B',\sigma')$ will now be rooted at $\Sigma\Delta$ and will have spine $\Theta\Delta$).

Finally, recall from sub-section \ref{rev-pair-def} that a tree pair $(A',B',\sigma')$ is a \emph{single rolling of type II
  from $(A,B,\sigma)$} if it is obtained from $(A, B,\sigma)$ by adding
a copy of a component $U$ of $A-B$ to $A$ at the last leaf in the
orbit of the repeller in $U$ and to $B$ at its image; or, by adding a
copy of a component $W$ of $B-A$ to $A$ at the first leaf in the
orbit of the attractor (the leaf of $A$ corresponding to the root node of $W$ in $B$) and to $B$ at its image.

The following lemma lists some basic properties of discrete train tracks.  The reader will not be required to use this lemma later in the paper, although it gives a separate view of some arguments.
 
 \bl
\label{conjugacyPreservesComponents}
Suppose $\tau_1 \sim P_1 =(D_1,R_1,\sigma_1)$ and $\tau_2 \sim P_2 = (D_2,R_2,\sigma_2)$ are elements of $V_n$, and that $TT_1$ and $TT_2$ are the corresponding train tracks derived from the revealing pair representatives $P_1$ and $P_2$ of these elements.  Then we have the following:
\begin{enumerate}
\item \label{conjugacyDTT} Suppose $f\in V_n$ so that $\tau_1^f = \tau_2$.  Then,
\begin{enumerate}
\item $f$ induces a $1-1$ correspondence between the components of $TT_1$ which describe dynamics around repelling/attracting orbits and the components of $TT_2$ which describe dynamics around repelling/attracting orbits.  
\item This correspondence also guarantees that the individual cycles in these components are also carried to cycles of the same length and type (repellers must move to repellers and attractors to attractors).
\item This correspondence preserves the labelings of the spine of the complementary components for corresponding cycles. 
\end{enumerate}
\item Suppose $\tau_1=\tau_2$ and $(D_1,R_1,\sigma_1)$ differs from $(D_2,R_2,\sigma_2)$ by an application of a rolling of type $X$.  Then
\begin{enumerate}
\item if $X$ is a rolling of Type E, then a train track connecting a source node to a sink node will be split into $n$ such parallel tracks along its length, or a periodic circle representing a periodic orbit of neutral leaves will be split into $n$ copies of ``parallel" periodic circles.
\item if $X$ is a single rolling of Type I, then there is a suffix $\Delta$ of the spine so that the labels of all the sub-arcs in the relevant repelling or attracting circle will be modified by the addition of $\Delta$ as a suffix to all of the labels on the circle.   Also, some source-sink chains which are incident on the affected circle will be lengthened by the length of the orbit of the repeller or the attractor.
\item if $X$ is a single rolling of Type II, then the source-sink chains with sources or sinks beginning or ending in the relevant complementary component will all increase their lengths by one (each of these new arcs will need appropriate labels added), and the sub-arc of the relevant repelling (attracting) circle to which the spine of the complementary tree is glued will move one location backward (forward) in the cyclic ordering of the arcs in that circle (respectively).  Finally, the labels of the affected arcs on the circle will change (the arc corresponding to the old spine will now be labelled by the old leaf label for the attractor or the repeller of that spine, while the arc corresponding to the new spine will have as root label its old label, and as leaf label, its old label concatenated with the word corresponding to the spine of the repeller of attractor).
\end{enumerate} 
\end{enumerate}  
\el
{\it Proof:} 
The latter points about rollings follow directly from the definitions of rollings, as given in subsection \ref{rev-pair-def}.

The first two sub-points of \ref{conjugacyDTT} are a result of the fact that $V_n$ is a group of homeomorphisms, and topological conjugacy preserves the properties mentioned.  The third sub-point of \ref{conjugacyDTT} follows from the fact that elements of $V_n$ do not change infinite suffixes; so, the points in the finite orbit of a repelling periodic point or of an attracting periodic point all have the same infinite repeating suffix (that is, as described in Lemmas \ref{AttractorProps} and \ref{RepellerProps}).  Now the conjugating element $f$ again cannot change this infinite suffix class, so the resulting orbit will consist of points with this same infinite suffix.  (Note first that that Belk and Matucci \cite{BelkMatucci}
also derive the infinite suffix of a finite repelling or attracting orbit as a conjugacy invariant for elements of $F<T<V$, using a method very similar to our discrete train tracks, although Belk and Matucci's definition
is mildly different, and second that one can choose a representative revealing pair for the conjugate version $\tau_2$ of $\tau_1$ so that the spine is cyclically rotated, but by applying a rolling of type $II$, this spine can be rotated back to the original spine, thus producing a discrete train track with the same spine labelling on the appropriate sub-arc of the circle (although, the prefixes of all nodes on the circle will be lengthened).)
\qquad $\diamond$

\subsection{Flow graphs}
The \emph{flow graph of any revealing pair $p = (A,B,\sigma)\sim \alpha\in
  V_n$} is a labeled directed graph which is, in some sense, a quotient object from a discrete train track (we loose the visual aspects of the branching of source sink lines as they leave the orbit of a repeller or join the orbit of an attractor, although even this information can be recovered from labels).  We now describe how to build a flow graph from a revealing pair.  
  
  For each repeller of $A$, we draw a vertex.  For each attractor of $B$, we draw a vertex.  For each
neutral leaf of $\labs$ that is part of a periodic orbit of neutral
leaves for $p$, we draw a vertex.  We draw a directed edge from a
repeller to an attractor for each source in the basin of repulsion of
the repeller whose iterated augmentation chain terminates in a sink in
the basin of attraction of the attractor (we call these
\emph{source-sink flow lines} or by similar language; they are in
one-one correspondence with the set of source-sink chains for the pair
$p$).  We draw a directed edge from each repelling and attracting
vertex to itself whenever the period of the corresponding repeller (or
attractor) is greater than one (we call this a \emph{repelling (or
  attracting) periodic orbit}).  We draw an edge connecting two
vertices representing periodic neutral leaves if a single iteration of
$\alpha$ will take the first leaf to the second, whenever these leaves are
not the same.  We label all source-sink flow lines with the
appropriate iterated augmentation chain.  We label all repelling and
attracting periodic orbits (even of length one) with the finite
periodic orbit of the actual points in $\CS_n$, each such point
labeled by its infinite descending path in $P_n$.  (Note that one can
detect the names of repelling and attracting basins by deleting the
infinite ``$\Gamma^{\infty}$''  portion of these labels.)

The diagram below is an example of a flow graph for the tree
$(A,B,\sigma)$ we have been examining. Strings of form $\Gamma^\infty$ 
are indicated by overlines in the labels of the diagram.

\begin{center}
\includegraphics[height = 100 pt, width = 400 pt]{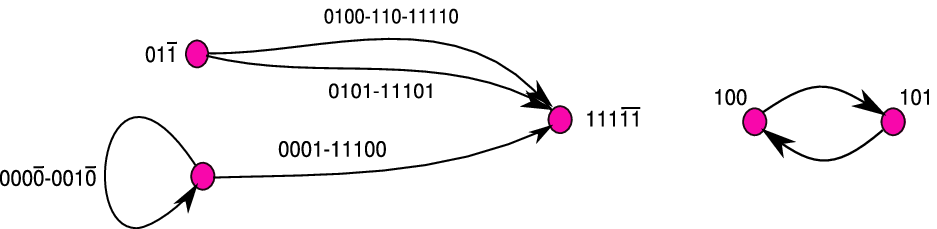}
\end{center}

The components of the flow graph for a revealing pair naturally
decompose into two sets; components representing the flow (under the
action of $\langle \alpha\rangle$ along the underlying sets of a
periodic orbit of neutral leaves for the revealing pair, called
\emph{torsion components}, and components representing flows from
basins of repulsion to basins of attraction, and characterizing the
orbits of repelling and attracting periodic points, called
\emph{non-torsion components}.

Any flow line in the graph from a repeller to an attractor can be
thought of as representing the complete bi-infinite forward and backward
orbit in $\CS_n$ of the interval underlying the source, under the action
of $\langle \alpha\rangle$.  The set so determined is called the
\emph{underlying set for the flow line}, and it limits on some of the
periodic repelling points of $\alpha$ (on one end), and on some of the
periodic attracting points of $\alpha$ (on the other end).  

Given a component $C$ of a flow graph for $\alpha$, we can discuss the
\emph{underlying set for the component $C$} in $\CS_n$.  We define
this set as the closure of the underlying sets of the flow lines of
the component if the component has a flow line (thus, we will capture the limiting repelling and attracting points in their finite orbits), and otherwise as the
underlying set of the finite periodic orbit of an appropriate neutral leaf.

We observe in passing the following remark, which the reader can
verify as a test of their understanding. The proof of this remark follows directly
from the definition of flow graph or can be obtained as a consequence of
Lemma \ref{conjugacyPreservesComponents}.

\begin{remark} \label{rollingAffects}

Given any revealing pair $p=(A,B,\sigma)$ representing $\alpha\in V_n$, and
$p'=(A',B',\sigma')$ the result of a single rolling of type II from $(A,B,\sigma)$,
the flow graph $\mathscr{F}_{\alpha,p'}$ for $\alpha$ generated by $p'$ is
identical to the flow graph $\mathscr{F}_{\alpha,p}$ for $\alpha$ generated by
$p$, except for one pair $(U,U')$ of corresponding components. The
components $U$ and $U'$ are both non-torsion components which have the
same underlying set and are isomorphic as graphs. However, the flow
lines of $U'$ will bear different labels from the flow lines of $U$, having had some of its iterated augmentation chains lengthened by one hop,
and, one of the periodic orbits in the labeling of $U'$, of either a
repelling periodic point or of an attracting periodic point, may be
cyclically permuted from the corresponding orbit as labeled in $U$.

\end{remark}

A mildly more difficult statement is that in the following lemma.  

\bl \label{flowsUnderPower}

Suppose $\alpha\in V_n$ is represented by a revealing pair $p
=(A,B,\sigma)\sim \alpha$.  If $z_1$ and $z_2$ underlie distinct non-torsion
components of $\mathscr{F}_{p,\alpha}$, then they cannot underlie the same
component of a flow graph representing $\alpha^k$ for any integer $k$.

\el 

{\it Proof:} 

Powers of $\alpha$ may split orbits of repelling (or attracting) periodic
points under the action of $\langle \alpha\rangle$, but they can never move
a point from the underlying set of one component of the flow graph to
the underlying set of other components; these sets are themselves
created as the unions of the images of the forward and backward orbits
of the points underlying the interval portions of the fundamental
domain of $\alpha$ which are carried in affine fashion by $\alpha$.

\qquad $\diamond$

The next lemma is a direct consequence of the work of Salazar in \cite{salazarPaper}.
\bl
\label{conjugationEffects}
Suppose $n$ is a positive integer and $\alpha\in V_n$ is represented by the
revealing pair $(A,B,\sigma)$.  
For all integers $k \ge 0$, set $m_k$ to be the 
number of finite periodic orbits of neutral leaves for the tree pair $(A,B,\sigma)$ with
period $k$ and set $r_k$ to be $m_k \pmod{n-1}$.

There is a conjugate $u = \alpha^f$ for
some $f\in V_n$, where $u$ has a representative tree pair $(E,F,\chi)$,
so that the following hold.
\be
\item The element $u$ has $r_k$ periodic orbits of
  neutral leaves with period $k$.
\item For each $k \ge 0$, the number $r_k$ is a conjugacy invariant of $\alpha$.
\item The number of components of $F-E$ is
the same as the number of components of $B-A$, and these components
have the same shape, and the same spine.  
\item Likewise, the number of
components of $E-F$ is the same as the number of components of $A-B$,
and these components have the same shape, and the same spine.
\ee
\el 
{\it Sketch of proof:}

The statements in the lemma above are for the most part stated
in Corollary 1 in Section 4.4 of \cite{salazarPaper}.  We explain our
statements relating to the modular arithmetic involved in the
reduction of the number of finite periodic orbits of the same length,
as that is not given Salazar's corollary.  The fact that spines are
preserved as well is not explicitly stated by Salazar, but it is an
immediate consequence of her techniques.

One can describe any conjugation in $V_n$ in terms of tree pairs by a
process which we now sketch.  Suppose $u = \alpha^f$ and let $f$ be
represented by a tree pair $(C,D,\tau)$.  There is an expansion
$(E',F',\chi')$ of $(A,B,\sigma)$ so that $C\subset (E'\cap
F')$.  One replaces the sub-tree $C$ inside of $(E'\cap F')$ by $D$,
producing a revealing tree pair $(E,F,\chi)$ for $u$ (one has to
remember to percolate out the effects of $\tau$ in both trees) with the same number
of neutral leaves as for the $(E',F',\chi')$ tree pair.

This process of conjugation produces a revealing pair for $\alpha^f$ whose
periodic neutral leaves are generated by expansions along the iterated
augmentation chains of the original neutral cycles (together with some
permuting of locations in the universal tree).  Thus, the number of
neutral periodic leaves of the result is in the same congruence class
modulo $(n-1)$ as for the tree pair for $\alpha$.

One can now find a conjugate of $\alpha$ which has all finite periodic orbits of neutral leaves of the same period adjacent in the universal tree $\mathscr{T}_n$. If there are more than $n$ such orbits, one can find a conjugate such that $n$ of the orbits travel in parallel, all carried in the $k$-orbit of an $n$-caret. One can also choose a conjugation so that one can build a new revealing pair for the resulting element with the same number of periodic orbit neutral leaves as in $(A,B, \sigma)$. By simple reductions along the full orbit of the caret, one can reduce the current representative tree pair to another revealing pair with $k \cdot (n-1)$ fewer periodic neutral leaves, and can repeat this overall process until there are fewer than $n$ orbits of neutral leaves for any period $k$.

\qquad$\diamond$

\section{A partition of $\CS_n$}
Given positive integer $n$ and some $\alpha\in V_n$ with
representative revealing pair $(A,B,\sigma)\sim\alpha$, we can
decompose $\CS_n$ as $T_{\alpha}\cup Z_{\alpha}$.  Here, we are using the following notation:
\begin{itemize}
\item $T_{\alpha}$: the subset of $\CS_n$ lying under the
neutral leaves in $A$ which are on cyclic orbits (possibly of length one).
\item $Z_{\alpha}$: the subset of $\CS_n$ underlying the root
  node of the complementary components of $A-B$ and $B-A$, as well
  as any neutral leaves which are part of source-sink chains.
\end{itemize}

The reader can easily verify the following.

\bl\label{domainDecomposition}
If $\alpha\in V_n$ and $\beta\in C_{V_n}(\alpha)$, then 
\be
\item $T_{\alpha}\beta = T_{\alpha}$, and
\item $Z_{\alpha}\beta = Z_{\alpha}$.
\ee
\el

This lemma allows us to work to comprehend centralizers over each set,
without regard to the behavior of these centralizers in the other
regions.  

From here to near the end of Section \ref{nontorsionSection}, we will assume that $n$ and
$\alpha$ are fixed, and $\alpha$ is the element whose centralizer in
$V_n$ we are analyzing, and that $(A,B,\sigma)$ is a revealing pair
representing $\alpha$.  We further assume through the use of the third
point of Remark \ref{conjugationRemark} and of Lemma
\ref{conjugationEffects} that for each periodic cycle length $m_i$ of
neutral leaves of $A$ under the action of $\langle \alpha\rangle$,
that there are precisely $r_i<n$ such cycles of periodic neutral
leaves.  These values $m_i$ and $r_i$ are the numbers which appear in
the semi-direct product terms in the left-hand direct product in the
statement of Theorem \ref{BigTheorem}.

Suppose $G\leq V_n$ and $X\subset \CS_{n}$.  We define
\[
G_X = \left\{v\in G \mid v|_{\CS_n\backslash X} = Id|_{\CS_n\backslash
  X}\right\}.
\]
so that $G_X$ is the subgroup of elements of $G$ which act as the
identity \textbf{except on} the set $X$. Lemma \ref{domainDecomposition}, assures
us that there are two commuting elements $\alpha_T\in {V_n}_{T_{\alpha}}$ and $\alpha_Z\in
{V_n}_{Z_{\alpha}}$ of $V_n$, so that $\alpha|_{T_\alpha} = \alpha_T|_{T_\alpha}$, and $\alpha|_{Z_{\alpha}} = \alpha_Z|_{Z_\alpha}$.  Thus, we see immediately that $\alpha = \alpha_T\alpha_Z$.

We will therefore restrict our attention to finding the centralizers
$C_{\,\,{V_n}_{T_\alpha}}(\alpha_T)$, and $C_{\,\,{V_n}_{Z_\alpha}}(\alpha_Z)$.  In
fact, we have the following corollary to Lemma
\ref{domainDecomposition}.

\bc \label{centralizerDecomp}

We have 
\[
C_{V_n}(\alpha) \cong C_{\,\,{V_n}_{T_\alpha}}(\alpha_T)\times C_{\,\,{V_n}_{Z_\alpha}}(\alpha_Z)
\] 

\ec

This explains the central direct product in our statement of Theorem \ref{BigTheorem}.

\section{Centralizers over the set $T_{\alpha}$\label{torsionCalc}}
Taking advantage of the decomposition given by Corollary \ref{centralizerDecomp}, we analyze the centralizer of $\alpha$ by restricting our attention to the set $T_\alpha$, over which $\alpha$ acts as an element of torsion.

Below, when we refer to a cycle of $\alpha$, we
mean a cycle of neutral leaves for $(A,B,\sigma)$.

\begin{lemma}
The leaves of $A$ over $T_{\alpha}$ can be partitioned into disjoint
sets according to cycle lengths, denoted by $S_{m_i}$, where each
$S_{m_i}$ consists of all the leaves of $A$ in a periodic orbit of
length $m_i$.  If $\beta\in C_{V_n}(\alpha)$, then $\beta$ preserves
the subsets of $\CS_n$ underlying the leaves in any particular set
$S_{m_i}$.  \end{lemma}

{\it Proof:} We first suppose $T_\alpha$ is not empty, and we further
suppose $\labs$ admits distinct neutral leaf cycles of length $k$ and
$m$ under the action of $\alpha$, where we chose our labels so that $k
< m$.  Finally, suppose $\beta\in C_{V_n}(\alpha)$ acts by mapping a
point $p'$ underlying a cycle of length $m$ to a point $p = p'\beta$
underlying a cycle of $\alpha$ of length $k$.  (If there is a
$\gamma\in C_{V_n}(\alpha)$ which maps a point underlying a cycle of
length $k$ to a point underlying a cycle of length $m$, then
$\gamma^{-1}$ will match our requirements.)  We now have the following
computation:

\[
p = p\alpha^k=p(\alpha^\beta)^k=p'\alpha^k\beta=q\beta.
\]
where $q=p'\alpha^k$ is not $p'$ since the orbit length for $p'$ under
$\alpha$ is $m$.  However, we have just shown that $q\beta = p$ and by
assumption $p'\beta = p$, so we have a contradiction.  \qquad$\diamond$

  We suppose throughout the remainder that there are $s$ distinct
  neutral leaf cycle lengths under the induced action of $\langle
  \alpha\rangle$ on the periodic neutral leaves in $\labs$, namely
  $\{m_1,m_2,\ldots,m_s\}$.

We thus can focus on how a particular $\beta\in C_{V_k}(\alpha)$ can
commute with $\alpha$ over the underlying set of the leaves in any
particular set $S_{m_i}$.  This is the reason for the left-hand direct
product with $s$ terms in our statement of Theorem \ref{BigTheorem}.
(Here, we are following the same logic as used in the beginning of
this section which allowed us to focus our attention on $T_\alpha$
based on the dynamical cause of the direct product decomposition of
Corollary \ref{centralizerDecomp}.)

For each leaf of $S_{m_i}$, we can consider its orbit in $S_{m_i}$
under the induced action of $\langle \alpha\rangle$.  Build a set
$F_{m_i}\subset S_{m_i}$ by taking one leaf from each such orbit.
Thus, $F_{m_i}$ is a collection of $r_i$ leaves.  Let $X_{m_i}$ denote
the subset of $\CS_n$ underlying $S_{m_i}$ and let $E_{m_i}$ denote
the subset of $\CS_n$ underlying the set $F_{m_i}$.  By construction
we see that a subset of the fundamental domain of $\alpha$ is
$D_{m_i}:=X_{m_i}/\langle\alpha\rangle\cong E_{m_i}$.

We now analyze the groups $G_{m_i}=C_{\,\,{V_n}_{X_{m_i}}}(\alpha)=
C_{V_n}(\alpha)\cap {V_n}_{X_{m_i}}$, which are individually
isomorphic to the terms in the left-hand direct product of Theorem
\ref{BigTheorem}.

As described in Section \ref{bigApproach}, $G_{m_i}$ is an extension
of its subgroup $K_{m_i}$ consisting of the elements in $G_{m_i}$
which have trivial induced action on $D_{m_i}$ and which act trivially
outside of $X_{m_i}$.

We now study $K_{m_i}$.  Let $\beta$ be an element of $K_{m_i}$, and
fix $\beta$ until we have finished our classification of $K_{m_i}$.

By the previous two paragraphs we see that $\beta$ must carry each set
underlying a leaf cycle in $S_{m_i}$ (under the action of
$\langle\alpha\rangle$) to itself.  Consider a leaf $\leafl\in
F_{m_i}$, and how $\beta$ moves points from the set underlying
$\leafl$ (fix this choice of leaf $\leafl$ for the remainder of the
discussion leading to the classification of $K_{m_i}$).  First, set
$\leafl_r = \leafl \alpha^r$, the $r$'th leaf in the orbit of $\leafl$
under the action of $\langle \alpha \rangle$ on $S_{m_i}$, for $r\in
\{0,1,\ldots, m_i-1\}$, and denote by $\tau_{\leafl,r}$ the set
underlying $\leafl_r$ in $\CS_n$.  If we take a point $p_0\in
\tau_{\leafl,0}$, its orbit under the action of $\langle
\alpha\rangle$ is $p_0$, $p_1 = p_0\alpha$, $p_2 = p_0\alpha^2$, and
etc., so that $p_r \in \tau_{\leafl,r}$.  In order for
$\langle\beta\rangle$ to have no induced action on $D_{m_i}$, we see
that $p_0\beta = p_c$ for some index $c$.  Now, in order to commute
with the action of $\langle\alpha\rangle$, we must have $p_r\beta =
p_{((r+c)\mod m_i)}$ for any index $0\leq r\le m_i$.
In particular,
$\beta$ must push the full orbit of $p_0$ under the action of $\langle
\alpha\rangle$ forward by some constant index less than $m_i$.  A
consequence of this is that any point $p \in E_{m_i}$
is itself in an orbit of
length less than or equal to $m_i$ under the action of $\langle
\beta\rangle$.  Extending this discussion as possible via recalling
our choices of $F_{m_i}$ and $\leafl$, we see that $\beta^{(m_i!)}$
must act trivially over the whole set $X_{m_i}$, so $\beta$ must be
torsion.

This now leads to the conclusion of our classification of $K_{m_i}$.
Suppose that $(C,D,\chi)$ is a revealing pair for $\beta$.  Since
$\beta$ is torsion, we see that $C=D$.  We assume (by taking a larger
revealing pair to represent $\beta$ if necessary), that $\leafl$ is a
node of $C$.  Since $\alpha$ takes the set $\tau_{\leafl,r}$ to the
set $\tau_{\leafl,((r+1)\mod {m_i})}$ in affine fashion, and $\beta$
commutes with $\alpha$ in such a way as to have no induced action on
$D_{m_i}$, it is straightforward to verify that the sub-tree
$T_{\leafl,r}$ in $C$ rooted at node $\leafl_r$ is identical in shape
to the sub-tree $T_{((r+1)\mod {m_i})}$ in $D$ rooted at
$\leafl_{((r+1)\mod {m_i})}$, for all indices $r$.  Further, by the
last sentence of the previous paragraph, for each leaf $\gamma_0$ of
$T_{\leafl,0}$ there is an integer $t_{\gamma_0} \ge 0$ such that
$\beta$ will send the corresponding leaf $\gamma_r$ of $T_{\leafl,r}$
to the corresponding leaf $\gamma_{((r+t_{\gamma_0})\mod m_i)}$ in
$T_{\leafl,((r+t_{\gamma_0})\mod m_i)}$.  
Any map $\CS_n \to \CS_n$ which is the identity outside of $X_{m_i}$ and which satisfies
these properties can be found in ${V_n}_{X_{m_i}}$, and a
straightforward topological argument (using the compactness of
$\CS_n$, and the basis cones of the topology on $\CS_n$. See Subsection \ref{TreesCantorSets}
for the definition) shows
that these maps form a subgroup which is isomorphic to the group
$Maps(\CS_n,\Z_{m_i})$. (That is, any continuous map from $\CS_n$ to
$\Z_{m_i}$ can be described as a rooted, finite, labeled $n$-ary
tree, where each label indicates where in $\Z_{m_i}$ to send the set
underlying the labeled leaf and this represents the offset in the orbit under $\langle 
\alpha\rangle$
for the element $\beta$ on the corresponding leaf.) 
Since the choice of map to $\Z_{m_i}$ on
the set $\tau_{\leafl,0}$ has no bearing (for the definition of
$\beta$) on the choice of map to $\Z_{m_i}$ for the sets underlying
the other leaves in $F_{m_i}$, we see that
\[
K_{m_i}\cong (Maps(\CS_n,\Z_{m_i}))^{r_i}.
\]

We now need to consider the structure of $G_{m_i}/K_{m_i} =
Q_{m_i}$. Thus, we are modding out the subgroup of elements of
${V_n}_{X_{m_i}}$ which commute with the action of
$\langle\alpha\rangle$ by the subgroup $K_{m_i}$.  In particular, we
are looking at the elements of ${V_n}_{X_{m_i}}$ which carry, for each
index $0 \le r \le m_i-1$, the sets underlying the $r$'th copy of the fundamental
domain $D_{m_i}$ to itself, where the map on the $r$'th copy of the
fundamental domain is precisely the conjugate version (under the
action of $\alpha^r$) of the map on the $0$'th copy $E_{m_i}$.
Therefore, the group $Q_{m_i}$ is
isomorphically represented by the restriction of the action of $V_n$
to the domain $E_{m_i}$.  In particular, $Q_{m_i} \cong
{V_n}_{E_{m_i}}$.  Since $E_{m_i}$ is given as the disjoint union of
$r_i$ distinct copies of $\CS_n$, we see that this group is precisely
the finitely presented group $G_{n,r_i}$ of Higman in
\cite{HigmanFPSG}!  That is, we have the following.

\[
Q_{m_i}\cong G_{n,r_i}
\]

Note that we can realize an isomorphic copy of $Q_{m_i}$ in $G_{m_i}$
as follows.  Let $\beta'$ be an element of ${V_n}_{E_{m_i}}$,
represented by the tree pair $(C',D',\chi')$.  The set $S_{m_i}$
decomposes as $m_i$ copies of the nodes of $F_{m_i}$ in the universal
tree $\ut{n}$ (take copy $C_r$ as the nodes in the set
$F_{m_i}\alpha^r$, using the induced action of $\langle \alpha
\rangle$ on the set of subsets of the nodes in $S_{m_i}$, for each
index $0\leq r < m_i$, and fix this definition of the sets $C_r$ for
the remainder of this subsection).  There is a revealing pair
$(C'',D'',\chi'')\sim\beta'$ expansion of $(C',D',\chi')$ which has
all of the nodes in the set $S_{m_i}$ as leaves of $C''$ and $D''$
(excepting the nodes in $F_{m_i}$, which themselves are roots of a
forest pair $f =(\mathscr{F}_d,\mathscr{F}_r,\theta)$ representing the
element of $G_{n,r_k}$ corresponding to $\beta'$).  By simply gluing a
copy of $f$ to each $C_r$ for $r>0$ (the nodes in these $C_r$ are
leaves of $C''$ and $D''$), we can build a new revealing tree pair
$(C,D,\chi)$ representing an element $\beta\in {V_n}_{X_{m_i}}$ which
acts on the set $\Gamma_r$ underlying $C_r$ as $\beta'$ acts on
$\Gamma_0$ under $C_0$, for each index $r$ (fix this definition of the
sets $\Gamma_r$ for the remainder of the section as well).  It is
immediate by construction that the group $\widehat{Q}_{m_i}$ of
elements $\beta$ so constructed is isomorphic with $Q_{m_i}$ and is a
subgroup of $G_{m_i}$ which splits the short exact sequence
$K_{m_i}\hookrightarrow G_{m_i}\twoheadrightarrow Q_{m_i}$. 

Thus, we have the following.

\[
G_{m_i}\cong K_{m_i}\rtimes Q_{m_i} \cong K_{m_i}\rtimes G_{n,r_i}
\]

We can complete our analysis of the centralizer of $\alpha$ over
$T_\alpha$ by showing the following lemma, which concludes the proof
of Corollary \ref{finiteGeneration} from the introduction.

\bl\label{finiteGenerationGmi}
The group $G_{m_i}$ is finitely generated.
\el

{\it Proof:}

First, we recall that for all positive integer values $n>1$ and $r$,
Higman's group $G_{n,r}$ is finitely presented (Theorem 4.6 of
\cite{HigmanFPSG}).  Let us denote by $\langle
A_{n,r}\,|\,R_{n,r}\rangle$ a finite presentation of $G_{n,r}$ for any
such $n$ and $r$.

We first describe our set of generators for $G_{m_i}$.  For each
generator $g'\in A_{n,r_i}$, we will take as a generator of $G_{m_i}$
the element $g\in \widehat{Q}_{m_i}$ which duplicates the effect of
$g'$ on $E_{m_i}=\Gamma_0$ over each of the sets $\Gamma_r$ underlying
the copies $C_r$ of the leaves $C_0$ over $\Gamma_0$, for each valid
index $r$.

At this stage, our collection of generators generates the group
$\widehat{Q}_{m_i}$, which is finitely presented still by carrying
over the relations of $G_{n,r_i}$ as well in corresponding fashion.
We need add only one further generator to generate the remainder of
$G_{m_i}$.  Let $g_1\in K_{m_i}$ be the element represented by the
revealing pair $(S,T,\theta)$.  We define $(S,T,\theta)$ as follows.
Let $\leafl$ be a node of $F_{m_i}$ (fix this choice and dependent
derived notation for the remainder of this subsection).  Let $S=T$ be
the minimal $n$-ary tree so that $\leafl_r = \leafl\alpha^r$ is a node
which is a parent of $n$-leaves of $S$, for all index values $0\leq r
< m_i$.  Let $\theta$ be the permutation from the leaves of $S$ to
the leaves of $T$ which takes the first child of $\leafl_r$ and sends
it to the first child of $\leafl_{((r+1)\mod m_i)}$, for all indices
$0\leq r<m_i$, and otherwise acts as the identity.  The element $g_1$
so constructed is our last generator.

We now show that $g_1$, together with our other generators, is
sufficient to generate $K_{m_i}$.  Each element of $K_{m_i}$
decomposes as a finite product of sub-node translations along the
$m_i$ orbit of $F_{m_i}$ in $S_{m_i}$ (under the action of $\langle
\alpha\rangle$).  That is, we choose a descendant node $p$ of
$F_{m_i}$ in the universal tree $\ut{n}$, and a translation constant
$0\leq t\le m_i$.  Then we translate the full orbit of $p$ under the
action of $\langle \alpha\rangle$ forward cyclically by the constant
$t$, while acting as the identity elsewhere.  We denote this
translation as $p_t$.  By choosing a specific $n$-ary forest rooted at
$F_{m_i}$ and translating each leaf of the forest in such a fashion,
we can obtain any element of $K_{m_i}$, as described above.  Now,
recall that $G_{n,r_i}$ acts transitively on the set of nodes in the
infinite $n$-ary forest descending from $F_{m_i}$ which do not happen
to represent the full domain of $G_{n_,r_i}$ (if $r_i = 1$, no element
of $V_n = G_{n,1}$ can take a proper sub-node to the root node).  Thus,
given any particular descendant node $p$ from a node of $F_{m_i}$ and
a translation distance $t$, we can find an element $\rho$ of
$\widehat{Q}_{m_i}$ taking the first descendant of $\leafl_0$ to $p$.
It is now immediate by construction that $p_1 = g_1^\rho$, and $p_t =
p_1^t$.  In particular, the set consisting of $g_1$ and the generators
of $\widehat{Q}_{m_i}$ together, is sufficient to generate $G_{m_i}$.
\qquad$\diamond$

\section{Centralizers over the set $Z_{\alpha}$\label{nontorsionSection}}

\newcommand{\rep}[1]{\mathcal{R}_{#1}}
\newcommand{\att}[1]{\mathcal{A}_{#1}}
\newcommand{\periodic}[1]{\mathcal{RA}_{#1}}

Let us fix a revealing pair $\mathfrak{p} = (A,B,\sigma)\sim \alpha$.
Let $\{\Gamma_1,\Gamma_2, \ldots,\Gamma_e\}$ represent the set of
non-torsion flow graph components of the flow graph
$\mathscr{F}_{\mathfrak{p},\alpha}$, where for each index $i$, we
denote by $X_i$ the component support of $\Gamma_i$.  Let $\alpha_i$
represent the element in ${V_n}_{X_i}$ such that $\alpha_i|_{X_i} =
\alpha|_{X_i}$, and suppose further that $\mathfrak{p}_i =
(A_i,B_i,\sigma_i)\sim\alpha_i$ is a revealing tree pair that is
identical to $\mathfrak{p}$ over the support $X_i$ of $\gamma_i$, so
that the flow graph $\mathscr{F}_{\mathfrak{p}_i,\alpha_i}$ is
identical to $\Gamma_i$.  Recall that by definition $Z_{\alpha}=\cup_i
X_i$. For any $\beta\in V_n$ set
\[
\periodic{\beta}:=\rep{\beta}\sqcup\att{\beta}.
\]
Note that for any such $\beta$, we have that $\periodic{\beta}$ is a finite
discrete set.

\begin{lemma}
\label{periodicAction}
Let $g \in C_{V_n}(\alpha)$.  
\begin{enumerate}
\item The group $\langle g\rangle$ acts on the set $\periodic{\alpha}$. 
\item Given $r\in \rep{\alpha}$, there is a basin of repulsion $U_r$ of $\alpha$ containing $r$ so that $U_rg$ is contained in a basin of repulsion $U_s$ for some repelling periodic point $s$ of $\alpha$.
\item Given $r\in \att{\alpha}$, there is a basin of attraction $U_r$ of $\alpha$ containing $r$ so that $U_rg$ is contained in a basin of attraction $U_s$ for some attracting periodic point $s$ of $\alpha$.
\item The group $\langle g\rangle$ acts bijectively on each of the sets $\rep{\alpha}$ and $\att{\alpha}$ of repelling and attracting periodic points of $\alpha$.
\end{enumerate}
\end{lemma}

{\it Proof:}
Let $p\in \periodic{\alpha}$ have an orbit of size $k$ under the action of $\langle \alpha\rangle$.  We now have
\[
pg = p\alpha^kg = pg\alpha^k
\]
Thus, $pg$ is fixed by $\alpha^k$, so that $\langle g\rangle$ bijectively
preserves $\periodic{\alpha}$.

The second point follows from the continuity of $g$ and the following
computation.  Let $r \in \rep{\alpha}$ so that $rg\neq r$ (if such a
repelling periodic point fails to exist, we automatically have the
second part of our lemma for the repelling periodic points). Choose
$\widetilde{U}_r \subseteq \CS_n$ an interval neighborhood of $r$
small enough so that $g$ is affine on $\widetilde{U}_r$ and so that
$\widetilde{U}_r$ is contained in a basin of repulsion for $\alpha$. Expand
the revealing tree pair $(A,B,\sigma)$ representing $\alpha$ by a
rolling of type II to create a new revealing tree pair $(A',B',\sigma')$
representing $\alpha$ with a complementary component $C$ rooted in some
node whose underlying set is contained in $\widetilde{U}_r$.  The root
node $U_r$ of $C$ represents a basin of repulsion for $\alpha$ which
is an interval neighborhood of $r$ carried affinely by $g$ to another
interval of $\CS_n$.  Assume that $r$ is in a periodic orbit of length
$k$ under the action of $\langle \alpha\rangle$.  Our result for repelling
periodic points follows easily from the following limit:
\[
(U_r)g\alpha^{-(nk)} = (U_r)\alpha^{-(nk)} g \to rg \textrm{\qquad (as }n\to \infty).
\]
Hence $g$ takes a basin of repulsion neighborhood of
$r$ into a neighborhood $N$ of some periodic repelling point $s = rg$ of
$\alpha$ with orbit length dividing $k$ where $N$ limits on $s$ under
powers $-(nk)$ of $\alpha$.  

A similar argument shows the third point of the lemma, and the final
point of the lemma is an immediate consequence of the previous three
points.  

\qquad$\diamond$

The following corollary is immediate from the first point of the lemma
above, together with the fact that $\alpha$ only admits finitely many
periodic points.  

\bc 

Let $g \in C_{{V_n}_{Z_\alpha}}(\alpha)$, then $\periodic{\alpha}\subset
\periodic{g}$.  

\ec

The following corollary now follows from Lemma \ref{periodicAction},
using the idea behind the proof from Section \ref{torsionCalc} that a
centralizer of an element of torsion must carry the set of
all finite orbits of length $k$ (under the action of the torsion
element) to itself.

\bc 
\label{pporbits}
Let $g \in C_{V_n}(\alpha)$, and let $r\in \rep{\alpha}$ or
$r\in\att{\alpha}$, with periodic orbit $(r_i = r\alpha^i)_{i =
  0}^{k-1}$ in $\rep{\alpha}$ or $\att{\alpha}$ respectively, then
there is a periodic orbit $(s_i)_{i = 0}^{k-1}$ in $\rep{\alpha}$ or
$\att{\alpha}$ respectively such that $s_i = r_ig$.

\ec

The next corollary depends on the proof of the second and third
points of Lemma \ref{periodicAction}.

\bc \label{commutingComponents}

Suppose $g\in
C_{V_n}(\alpha)$ sends some point $z\in X_i$ to a point
$zg\in X_j$ for some indices $i$ and $j$, then
$g$ will send $X_i$ bijectively to $X_j$.

\ec {\it Proof:} 

We first show that $g$ will send all of the periodic repelling and
attracting points $\alpha$ within $X_i$ into $X_j$.  

Set $P_i:= X_i\cap\periodic{\alpha}$.

If $g$ takes two points in $P_i$ to the supports of distinct
components of the flow graph of $\alpha$, then $X_i$ must admit a flow
line which has some point $r_1$ in its periodic repelling orbit vertex
label sent to a periodic repelling point $r_1g$ of one component
$\Gamma_k$ while having another point $a_1$ in its periodic attracting
orbit vertex label sent to a periodic attracting point for $\alpha$ in
a label of a distinct component $\gamma_m$.  But now, there are basins
of repulsion $U_1$ and attraction $W_1$ around $r_1$ and $a_1$
respectively which are carried by $g$ into $X_k$ and $X_m$
respectively.  This last is a contradiction, as follows.

Recall that some non-zero power $k$ of $\alpha$ fixes all of the repelling
and attracting periodic points of $\alpha$.  Take $p_1 \in
U_1\backslash\{r_1\}$.  It must be that for all integers $z$,
\[
p_1g\alpha^{kz} = p_1\alpha^{kz}g,
\]
however, for $z$ large and negative, $p_1\alpha^{kz}$ is near $r_1$ while
for $z$ large and positive, $p_1\alpha^{kz}$ is near $a_1$.  In particular
$\alpha^k$ has a flow line connecting $r_1g$ to $a_1g$, which is not
possible by Lemma \ref{flowsUnderPower}.

If $p_1$ is in the support of a flow line, then by considering powers
of $\alpha$ using similar arguments as above we can show the whole flow
line is sent by $g$ to a single flow graph component of $\alpha$.

\qquad $\diamond$

We now define the function
\[
\mathcal{S}(g) = \log_{2n-1} \left( \prod_{r \in \rep{\alpha}} rg' \right).
\]
where $rg'$ denotes the slope of $g$ at the the repeller $r$.
This map is well-defined by using the recognition that $g$ is affine
in small neighborhoods of points in $Per(\alpha)$.  It is not too hard to
see that the function $\mathcal{S}: V_n \to \mathbb{Z}$ is not a
homomorphism in general.

\begin{lemma}
\label{SGroupHom}
The map $\mathcal{S}: C_{V_n}(\alpha) \to \mathbb{Z}$
is a group homomorphism.
\end{lemma}

{\it Proof:} Given $g_1,g_2 \in C_{V_n}(\alpha)$ we compute
$\mathcal{S}(g_1 g_2)$ directly from the definition.  We note that in
small interval neighborhoods of the points in $\rep{\alpha}$, the maps
$g_i$ are differentiable, and so we can apply the chain rule. Since
$g_i$ acts bijectively on the set $\rep{\alpha}$ (by Corollary \ref{pporbits})
each of the two terms of the product appears exactly once and so
$\mathcal{S}(g_1 g_2) = \mathcal{S}(g_1) +\mathcal{S}(g_2)$
\qquad$\diamond$

\bigskip
We will now shift attention to the local behavior of $C_{{V_n}_{X_i}}$
over the region $X_i$ for any particular index $i\in\{1,2,\ldots,t\}$.
For each such index $i$, set
\[
\mathcal{S}_i(g) = \log_{2n-1} \left( \prod_{r \in \rep{\alpha_i}} rg' \right).
\]
By the previous lemma, $\mathcal{S}_i$ is a group homomorphism when we
restrict its domain to either $C_{{V_n}_{X_i}}(\alpha_i)$ or even to
$C_{V_n}(\alpha)$, and further, if $\beta_i\in
C_{{V_n}_{X_i}}(\alpha_i)$ while $\beta\in C_{V_n}(\alpha)$, with
$\beta_i|_{X_i} = \beta|_{X_i}$, then $S_i(\beta_i) = S_i(\beta)$.
For our immediate purposes below, we will use
$S_i:C_{{V_n}_{X_i}}(\alpha_i)\to\mathbb{Z}$.

The previous lemma now immediately implies the
existence of the following exact sequence:
\begin{equation}\label{eq:second-exact}
0 \to \ker(\mathcal{S}_i) \hookrightarrow C_{\,\,{V_n}_{X_i}}(\alpha_i) \twoheadrightarrow \mathrm{im}(\mathcal{S}_i)=\mathbb{Z} \to 0
\end{equation}

\begin{lemma}(Stair Algorithm)
\label{staircase}
Let $g_1,g_2
\in C_{{V_n}_{X_i}}(\alpha_i)$ and $r,s$ be periodic repelling points of $\alpha_i$, for some index $i$. If
\[
rg_1=rg_2 = s
\]
and the slope of $g_1$ at $r$ is equal to the slope of $g_2$ at $r$ , then
$g_1=g_2$.
\end{lemma}

{\it Proof:} 

There is a basin of repulsion $U_r$ of $\alpha_i$ containing $r$ so that $g_1=g_2$ on $U_r$.  
Let $x\in U_r\backslash\{r\}$.  Since
\[
xg_1\alpha_i^n = xg_2\alpha_i^n = x\alpha_i^ng_1 = x\alpha_i^ng_2
\] 
for all integers $n$, we see that on the underlying
support of any flow line $L$ of $\alpha_i$ limiting on $r$ we have $g_1|_L =
g_2|_L$.  But this now means that $g_1$ agrees with $g_2$ on all of the
underlying support of the edges of the graph $\Gamma_i$.  In particular,
$g_1=g_2$ on sets limiting to each of the attracting orbits on the other
ends of the flows lines of $\Gamma_i$ leading away from the repelling orbit of $r$.
Since $g_1$ and $g_2$ are always affine in small neighborhoods of the
repelling and attracting periodic points for $\alpha$, we then see that $g_1$
and $g_2$ actually agree on small neighborhoods of the
attracting periodic orbits which appear on the terminal ends of those
flow lines leaving the orbit containing $r$.

We now repeat this argument moving away from the attracting orbits to
new repelling orbits along new flow lines, where we again have that
$g_1 = g_2$ along these flow lines.  Now by the connectivity of
$\Gamma_i$, $g_1=g_2$ over $X_i$.

\qquad$\diamond$

\begin{lemma}
There is an exact sequence
\begin{equation}\label{eq:first-exact}
0 \to \mathbb{Z}{\to} C_{\,\,{V_n}_{X_i}}(\alpha_i) \overset{q}{\to} Q \to 0
\end{equation}
where $Q \le \mathrm{Sym}(\mathcal{R}_{\alpha_i})$. In particular, $C_{\,\,{V_n}_{X_i}}(\alpha_i)$ is virtually infinite cyclic.
\end{lemma}

{\it Proof:}
Define $M:=\{\beta \in C_{\,\,{V_n}_{X_i}}(\alpha_i) \mid \beta(r)=r, \, \, \forall r \in \mathcal{R}_{\alpha_i}\}$.
Fix $r_1 \in \mathcal{R}_{\alpha_i}$ and define the following map
\[
\begin{array}{cccc}
\varphi: & M & \longrightarrow & \mathbb{Z} \\
         & \beta & \longmapsto & \log_{2n-1} (r_1\beta').
\end{array}
\]
By Lemma \ref{staircase}, the map $\varphi$ is injective and so $M
\cong \mathbb{Z}$. Now we observe that $M$ is the kernel of the action
of $C_{\,\,{V_n}_{X_i}}(\alpha_i)$ on
$\mathcal{R}_{\alpha_i}$ and so we get a natural map
$q:C_{\,\,{V_n}_{X_i}}(\alpha_i)\to Q$ where $Q:=
C_{\,\,{V_n}_{X_i}}(\alpha_i) /M \le
\mathrm{Sym}(\mathcal{R}_{\alpha_i})$.
\qquad$\diamond$

We recall a few elementary facts from group theory. The proofs are
immediate, but we provide them for completeness.

\begin{lemma}
\label{omegaLemma}
Let $G,H$ be groups such that the following sequence
\[
0 \to H \to G \overset{\omega}{\to} \mathbb{Z} \to 0.
\]
is exact. Then $G \cong H \rtimes \mathbb{Z}$.
\label{ZSplit}
\end{lemma}
{\it Proof:} Let $z \in \omega^{-1}(1)$ then, by construction,
$\mathbb{Z} \cong \langle z \rangle \le G$ and $H \cap \langle z \rangle =
\{1_G\}$. If $g \in G$ and $\omega(g)=k$, then $z^k g^{-1} \in H$, and
so $G = H \langle z \rangle$. \qquad $\diamond$

\bigskip
By Lemmas \ref{SGroupHom} and \ref{ZSplit} applied on the exact
sequence (\ref{eq:second-exact}) we obtain that
$C_{\,\,{V_n}_{X_i}}(\alpha_i) \cong \ker(\mathcal{S}_i) \rtimes
\mathbb{Z}$.  We will now show that $\ker(\mathcal{S}_i)$ is a finite
group and that it coincides with the set of torsion elements of
$C_{\,\,{V_n}_{X_i}}(\alpha_i)$.

\begin{lemma} 

Let $G,Q,K,C$ be groups, where $C=\mathbb{Z}$ and $Q$ is a finite
group of order $m$. Assume the following two sequences
\[
0 \to C=\mathbb{Z} \overset{\varphi}{\to} G \overset{\psi}{\to} Q \to
0
\]
and
\[
0 \to K \overset{\tau}{\to} G \overset{\sigma}{\to} \mathbb{Z} \to 0
\]
are exact. Then $K$ is a finite group and $G \cong K \rtimes \mathbb{Z}$.
\label{two-exact-sequences}
\end{lemma}

{\it Proof:} 

In this proof we write $M \le_f N$ to denote that $M$ is a finite
index subgroup of a group $N$.  Let $h \in \sigma^{-1}(1)$ and
$I:=\langle h \rangle \le G$. By Lemma \ref{ZSplit} we have $G \cong K
\rtimes I \cong K \rtimes \mathbb{Z}$. We need to show that $K$ is a
finite group.  By assumption, $C \le_f G$ and so we observe that
\[
\frac{I}{C\cap I} \cong \frac{IC}{C} \le \frac{G}{C} \cong Q,
\]
therefore implying that $C \cap I \le_f I$. In particular, $C \cap I$
is a non-trivial group, hence $C\cap I \le_f C$ is too. By definition,
for every $g \in G$, we have $g^m \in C$.  Since $C\cap I \le_f C$,
there is an integer $k$ such that $g^{km} \in C \cap I$, for every $g
\in G$.

If $g \in K$, we have that $g^{km} \in K \cap C \cap I \le K\cap I=
0$. Therefore, $K$ is a torsion subgroup of finite exponent, hence $K
\cap C = 0$ and so the first exact sequence implies that $K \cong
\psi(K) \le Q$ and therefore it is finite. 

\qquad$\diamond$

\bigskip

Applying Lemma \ref{two-exact-sequences} on the two exact sequences
(\ref{eq:second-exact}) and (\ref{eq:first-exact}) we deduce the
following result:

\begin{corollary}
Let $i\in \{1,2,\ldots, e\}$.  The centralizer
$C_{\,\,{V_n}_{X_i}}(\alpha_i)$ is isomorphic to $\ker(\mathcal{S}_i)
\rtimes \mathbb{Z}$ and $\ker(\mathcal{S}_i)$ is a finite group.
\label{thm:kernel-finite-group}
\end{corollary}

The following is a fairly technical statement about roots which does
not particularly assist us in our exploration of the centralizer of
$\alpha$, but which deserves to be stated.  It represents a viewpoint
on the underlying reason for the corollary on page 68 of
\cite{HigmanFPSG}.  Of course, a separate argument can also be given
simply by noting how the lengths of spines of complementary components
for representative revealing pairs change when one takes an element to powers.

\begin{corollary}
Suppose $\beta\in V_n$ is a non-torsion element with precisely one flow graph component $\Upsilon$ which has component support $Y$.  Then the set of roots of $\beta$ in ${V_n}_Y$ is finite.
\end{corollary}
{\it Proof:} Consider the exact sequence from Lemma \ref{omegaLemma},
applied to $C_{{V_n}_Y}(\beta)$.  If $\beta\in \omega^{-1}(k)$ for
some positive integer $k$, then all of the roots of $\beta$ in
${V_n}_Y$ occur in the sets $\omega^{-1}(j)$ for integers $j$ which
divide $k$.  Thus the roots all occur in a finite collection of finite
sets.  \qquad$\diamond$

We note in passing that it may be the case that for all $\gamma\in
\omega^{-1}(1)$ we have that $\gamma^k = \tau_\gamma\cdot \beta$,
where $\tau_\gamma$ is a non-trivial torsion element for each such
$\gamma$.

We now analyze the kernel of our $\mathcal{S}_i$ homomorphism a bit further.

\begin{lemma}
Let $i\in \{1,2,\ldots,e\}$.  The kernel $\ker(\mathcal{S}_i)=\{g \in
C_{{V_n}_{X_i}}(\alpha_i) \mid g^k=id, \; \; \mbox{for some} \; \; k \in
\mathbb{Z}\}$, that is, the kernel of $\mathcal{S}_i$ is the set of
torsion elements in $C_{{V_n}_{X_i}}(\alpha_i)$.
\end{lemma}

{\it Proof:} Let $T$ be the set of torsion elements in
$C_{{V_n}_{X_i}}(\alpha_i)$.  A priori, this may not be a subgroup.

We observe that $T \subseteq \ker(\mathcal{S}_i)$, because the slope of
$g \in T$ multiplies to one across the full cycle of every periodic
orbit of neutral leaves of $g$, and the orbits of the repelling
periodic orbits of $\alpha_i$ are carried to each other by the action of any
$\langle g\rangle$ for any $g\in C_{{V_n}_{X_i}}(\alpha_i)$.  So, for $g\in T$, we
have by the definition of $\mathcal{S}_i$ that $\mathcal{S}_i(g)=0$.

By the previous Corollary, $\ker(\mathcal{S})$ is a finite group,
hence $\ker(\mathcal{S}) \subseteq T$.  \qquad$\diamond$

\bigskip

At this juncture, we have pushed our analysis of the centralizer of an element of $V_n$ with a discrete train track (or flow graph) with one component, which represents a non-torsion component, as far as necessary for us to be able to support the structure described in the right hand product of our Theorem \ref{BigTheorem}.  

{\it Proof of the right hand product structure of Theorem \ref{BigTheorem}:}

We now partition the non-torsion flow graph components $\{\Gamma_i\}$ by the rule
$\Gamma_i \sim \Gamma_j$ if there exists $f\in C_{V_n}(\alpha)$ such that $(X_i)f=X_j$. In this case we note
that $\alpha_i^f=\alpha_j$ and therefore $C_{V_{X_i}}(\alpha_i) \cong C_{V_{X_j}}(\alpha_j)$.
In particular, we have that $\ker(\mathcal{S}_i) \cong\ker(\mathcal{S}_j)$.

Let us call such supports ($X_*$ which can be carried to each other by an $f \in C_{V_n}(\alpha)$) \emph{supports of isomorphic connected components}.  Let $\{ICC_1, ICC_2, \ldots, ICC_t\}$ be the set of $\sim$-equivalence classes
of isomorphic (flow graph) connected components (each of which is denoted by
$ICC_j$, for some $j$), where $q_j$ is the cardinality of the set $ICC_j$. 
For each $j$, order the members of $ICC_j$ and let $\Gamma_{i,j}$ be the $i$-th member of $ICC_j$. We re-index the $X_i, \, \alpha_i, \, \mathcal{S}_i$ using the new double-index 
notation in corresponding fashion.  Let us fix conjugators $\gamma_{1,k,j}$ in $C_{V_n}(\alpha)$ carrying $X_{1,j}$ to $X_{k,j}$ (acting as the identity elsewhere), and use these to generate a full permutation group $P_{q_j}$ which acts on the underlying sets $X_{*,j}$ of the elements in $ICC_j$. 
Set $A_j := \ker(\mathcal{S}_{1,j})$ and note that $\ker(\mathcal{S}_{i,j})=A_j^{\gamma_{1,i,j}}$.

For fixed $j$, we then have the centralizer of $\alpha$ over the set $\cup_i X_{i,j}$ consists of any self centralization on each individual $X_{i,j}$ (this group will be congruent
to $B_{i,j}:=\ker(\mathcal{S}_{i,j}) \rtimes \mathbb{Z}$), together with a product by any permutation of these components.  Thus, this group is 
\[
G_j \cong \left(\prod_{i=1}^{q_j} B_{i,j} \right) \rtimes  P_{q_j} = \left(A_j \rtimes \mathbb{Z} \right) \wr P_{q_j}.
\]
(note that the product $\prod_{i=1}^{q_j} B_{i,j}$ is normal in $G_j$).
Now as non-isomorphic non-torsion components cannot be mapped onto each other by the action of $f$, we see that the action of $f$ on the non-torsion components of $\alpha$ is describable as an element from the direct product of the individual groups $G_j$.  Since $t$ is
the number of isomorphism classes $\{ICC_j\}$,
we obtain our statement of Theorem \ref{BigTheorem}.

\qquad $\diamond$

Below is a revealing pair $(A,B,\sigma)$ representing element $\alpha$ of
$V_2$, which has centralizer congruent to $(\Z_2\times\Z_2)\times\Z$;
the group corresponding to $A_1$ is the Klein $4$-group.  This example
is included to answer a question of Nathan Barker which arose in
conversation.

\begin{center}
\psfrag{1}[c]{$\,\,1$}
\psfrag{2}[c]{$\,\,2$}
\psfrag{3}[c]{$\,\,3$}
\psfrag{4}[c]{$\,\,4$}
\psfrag{5}[c]{$\,\,5$}
\psfrag{6}[c]{$\,\,6$}
\psfrag{7}[c]{$\,\,7$}
\psfrag{8}[c]{$\,\,8$}
\psfrag{9}[c]{$\,\,9$}
\psfrag{10}[c]{$\,\,10$}
\psfrag{11}[c]{$\,\,11$}
\psfrag{12}[c]{$\,\,12$}
\psfrag{13}[c]{$\,\,13$}
\psfrag{14}[c]{$\,\,14$}
\psfrag{15}[c]{$\,\,15$}
\psfrag{16}[c]{$\,\,16$}
\psfrag{17}[c]{$\,\,17$}
\psfrag{18}[c]{$\,\,18$}
\psfrag{19}[c]{$\,\,19$}
\psfrag{20}[c]{$\,\,20$}
\psfrag{\alpha}[c]{$\alpha$}
\includegraphics[height=150 pt,width=300 pt]{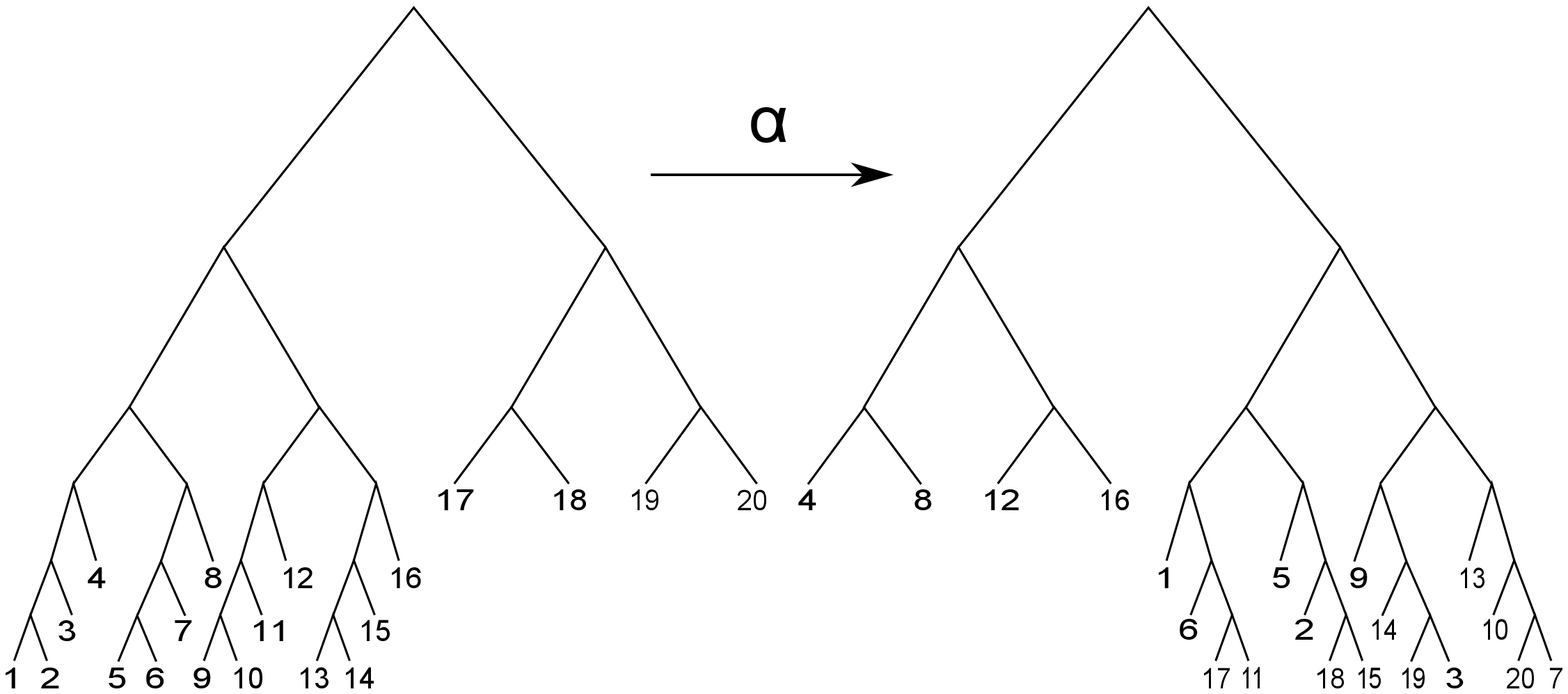}
\end{center}

Below, we sketch the flow graph associated with $(A,B,\sigma)$.  We do
not label the flow lines with their iterated augmentation chains, as
that will only serve to clutter the essential aspects of this graph.
The three different types of flow lines for this revealing pair are
encoded by the different methods used in drawing the directed edges
(solid, versus two distinct flavors of dashing).

\begin{center}
\psfrag{0}[c]{$Y$}
\psfrag{1}[c]{$X$}
\includegraphics[height=150 pt,width=300 pt]{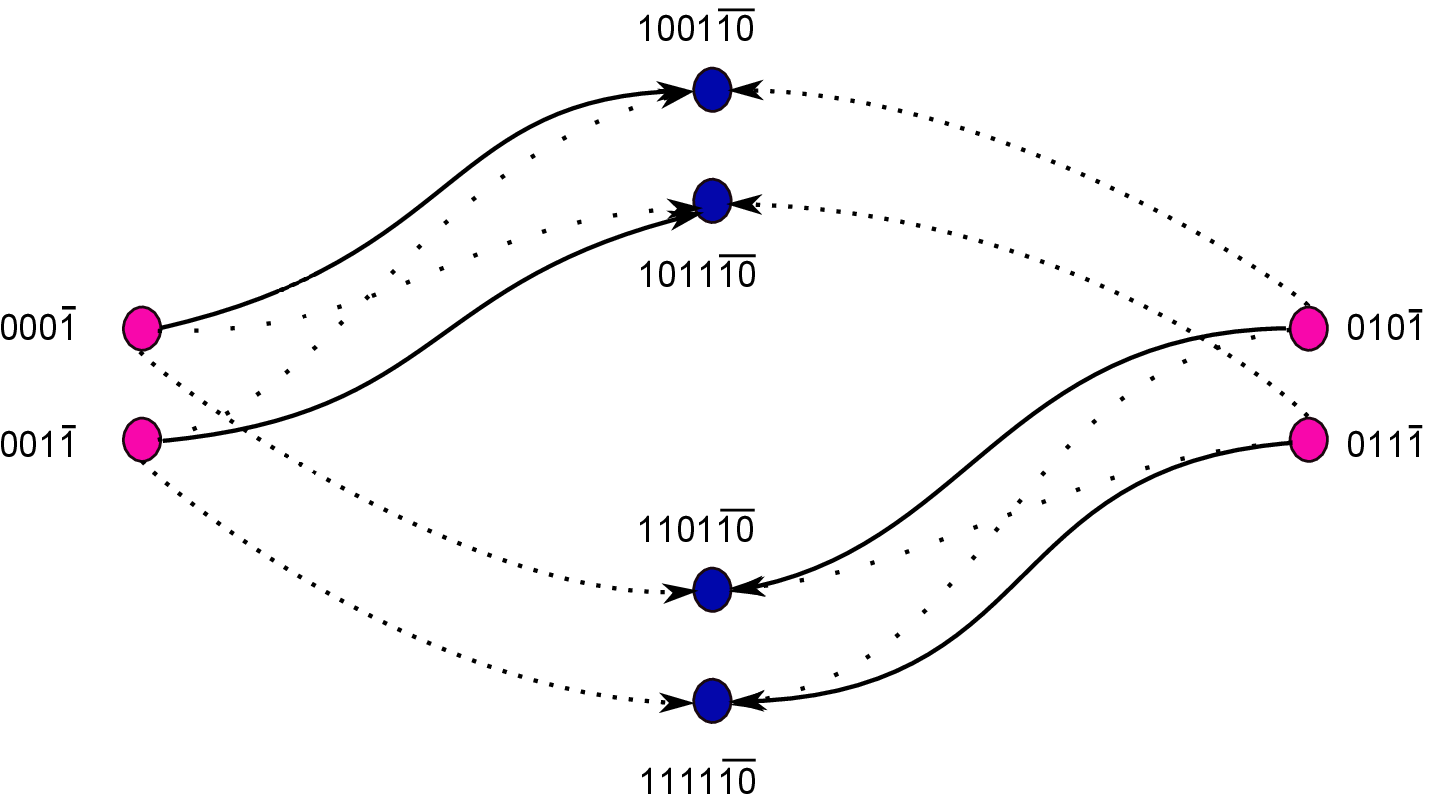}
\end{center}

In this graph, the repelling fixed points label the four vertices to
the left and right, while the attracting fixed points are the four
vertices in the middle.

Graph connectivity and the constraint on preserving flow-line types
provide enough information to guarantee that if $\tau$ is a
non-trivial torsion element commuting with $\alpha$, then $\tau$ is
completely determined by which repelling fixed point one chooses to
send $000\overline{1}$ to.  All such elements are order two, and the
elements $\beta$ and $\gamma$ below generate the Klein $4$-group $K$
so determined by the orbit dynamics on the repelling fixed points.
Thus, $A_1\leq K$.

\begin{center}
\psfrag{1}[c]{$\,\,1$}
\psfrag{2}[c]{$\,\,2$}
\psfrag{3}[c]{$\,\,3$}
\psfrag{4}[c]{$\,\,4$}
\psfrag{\beta}[c]{$\beta$}
\includegraphics[height=100 pt,width=200 pt]{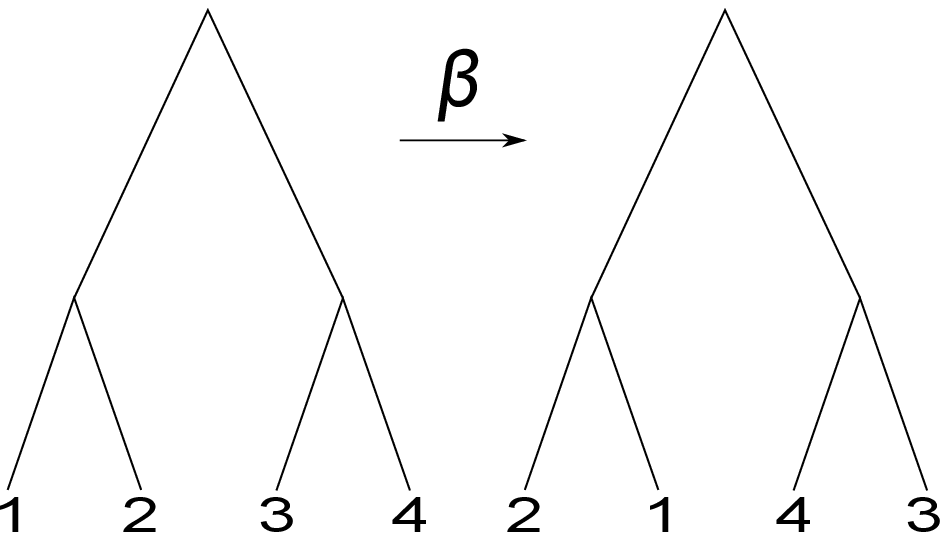}
\qquad\qquad\qquad
\psfrag{1}[c]{$\,\,1$}
\psfrag{2}[c]{$\,\,2$}
\psfrag{3}[c]{$\,\,3$}
\psfrag{4}[c]{$\,\,4$}
\psfrag{5}[c]{$\,\,5$}
\psfrag{6}[c]{$\,\,6$}
\psfrag{7}[c]{$\,\,7$}
\psfrag{8}[c]{$\,\,8$}
\includegraphics[height=100 pt,width=200 pt]{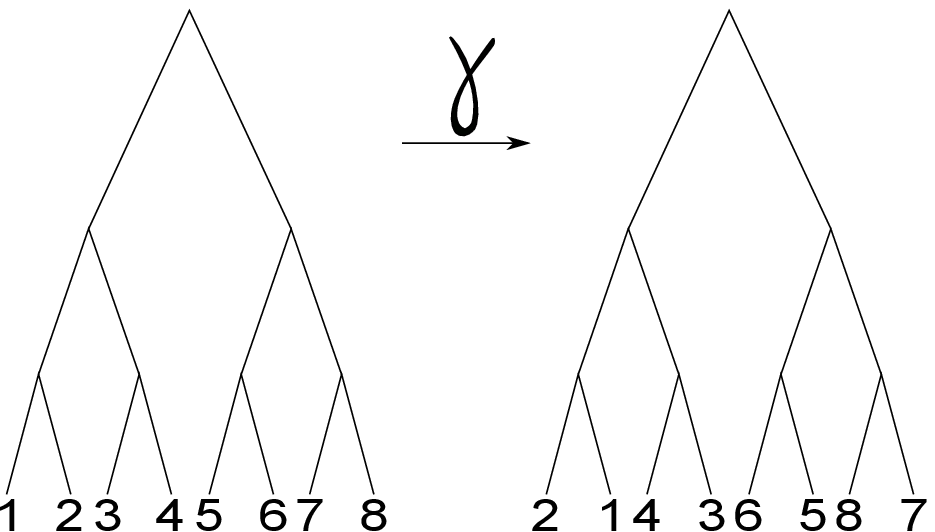}
\end{center}

Since $\beta$ and $\gamma$ can be realized as a group of permutations
of the leaves of the common tree $A\cap B$ which are in the repeller
and attractor iterated augmentation chains, while preserving the orbit
dynamics of $\langle \alpha\rangle$ over $\periodic{\alpha}$, we see that
every element in $K$ commutes with $\alpha$.
Note that as $\alpha$ has no roots, the generator of the $\mathbb{Z}$ factor
is $\alpha$, so $A_1 \rtimes \mathbb{Z} \cong A_1 \times \mathbb{Z}$, in this example.

\section{Cyclic subgroup (non)distortion}

We assume the reader is familiar with distortion of subgroups in a
group.  We use definitions consistent with their usage in
\cite{FarbDistortion}.  

The calculation of the (non)distortion of the cyclic subgroups in
$V_n$ is very close is spirit to the calculations in \cite{BurilloDistortion,
cWladis1, BCSTCombinatoricsT}.  However, by using properties of revealing
pairs, we remove many of the technical obstructions usually associated with
counting carets and this shortens the arguments in those papers. 

We believe that technology such as the revealing pair technology should generally simplify
proofs of non-distortion of cyclic subgroups in the various families of generalized
Thompson's groups following arguments similar to that which is given below.
It would be interesting to see a general such tool developed for 
these families of groups.

We are indebted to M. Kassabov for a discussion of the spirit
of this type of argument.

{\it Proof of Theorem \ref{distortion}:}

Suppose $\alpha\in V_n$ and that $\langle \alpha \rangle\cong \Z$.  Suppose
further that $X = \{x_1, x_2, \ldots, x_q\}$ is a finite generating
set of $V_n$ which is closed under inverses.

There is a minimal positive integer $P_X$ so that for all $p\in
\CS_n$, the slope of each $x_j$ in small neighborhoods of $p$ is
$(2n-1)^{s_p}$ where $|s_p|<P_X$.  

Suppose $(A,B,\sigma)\sim \alpha$, and let
$r$ be a repeller for this tree pair.  Suppose further the
iterated augmentation chain of $r$ is ${r} = {r}_0$,
${r}_i = {r_0}{\alpha}^i$ for $0\leq i\leq u$ (so that
${r}_u$ is the root of the complementary component ${C_
r}$
of ${A} - {B}$ containing the repeller
${r}$).  Let ${\Gamma_r}$ be the spine of ${C_r}$, and
suppose the length of ${\Gamma_r}$ is ${L_r}$.

For each index $0<i<u$, there is a jump $J_{(r,i)}$ in depth in the infinite
binary tree, from the depth of ${r}_i$ to ${r}_{i+1}$. Further, set
$J_{(r,0)}$ to be the jump from the depth of ${r}_u$ to the depth of
${r}_1$. (E.g., if ${r}_u$ has depth $4$, and ${r}_{1}$ has
depth $3$, then we set $J_{(r,0)}$ to be $1$.)  Given a positive integer
$z$, set 
\[
S_{(r,z)} = \Sigma_{e = 1}^{(z\!\!\!\!\mod u)}J_{(r,(e-1))},
\]
 the partial sum of the first $z\!\!\mod u$ jumps. (This sum may be
 negative.)  Note that the sum of all $u$ jumps is zero, so that
 
\[
\Sigma_{e = 1}^z J_{(r,\,(e-1) \!\!\!\!\mod u)} = S_{(r,z)}.
\]

Now fix a particular positive integer $z$.  Set $w = \lfloor
z/u\rfloor$, the largest positive integer less than or equal to $z/u$.
Now, if ${y}_i$ is the repelling periodic point under the leaf ${r}_i$
for $0\leq i<u$, direct calculation shows that the slope of $\alpha^z$ at
${y}_0$ is
\[
((2n-1)^{{L_r}})^{w+1}\cdot(2n-1)^{S_{(r,z)}}.
\]

By the chain rule, we require at least $\lceil ((L_r\cdot
w+L_r+S_{(r,z)})/P_X\rceil$ generators from $X$ to create the element
${\alpha}^z$ (where $\lceil \cdot \rceil$ denotes
the largest positive integer greater than or equal to 
its argument).  This last function can be interpreted as a function in $z$
that is bounded below by an affine function $g:\N\to\Q$ where $g(z)=
(1/W)\cdot z + O$ with positive slope $1/W<(L_r/(u\cdot P_X))$ for
some integer $W$ and vertical offset $O\in \Z$ (for technical reasons, 
choose $O$ so that $g(0)<0$).  Now the function $f_{pos} = g^{-1}$
(restricted and co-restricted to $\N$) is an affine distortion
function for the positive powers of $\alpha$ in $\langle \alpha \rangle$ within
$V_n$ ($\alpha^z$ has minimal length $z$ when expressed as an element in
$\langle \alpha\rangle$ using the generating set $\{\alpha,\alpha^{-1}\}$, and by
construction, $z< f_{pos}(m)$, where $m$ is the minimal word length of
$\alpha^z$ as expressed in $X$).

A similar argument produces an affine distortion function
$f_{neg}:\N\to\N$ for the negative powers of $\alpha$.  Thus,
$f:=f_{pos}+f_{neg}$ will be an affine distortion function for the
whole of $\langle \alpha\rangle$ in $V_n$.  In particular, $\langle
\alpha\rangle$ is undistorted in $V_n$.

\qquad$\diamond$ 

\bibliographystyle{amsplain}
\def\cprime{$'$}
\providecommand{\bysame}{\leavevmode\hbox to3em{\hrulefill}\thinspace}
\providecommand{\MR}{\relax\ifhmode\unskip\space\fi MR }
\providecommand{\MRhref}[2]{%
  \href{http://www.ams.org/mathscinet-getitem?mr=#1}{#2}
}
\providecommand{\href}[2]{#2}

\end{document}